\newif\iffinal
\finaltrue	

\documentclass[letterpaper,12pt,reqno]{amsart}
\RequirePackage[utf8]{inputenc}
\usepackage[portrait,margin=2.6cm]{geometry}
\iffinal\else\usepackage[notref,notcite]{showkeys}\fi
\usepackage{mathrsfs,soul}
\usepackage{hyperref}
\usepackage[foot]{amsaddr}
\usepackage{amssymb,amsthm,amsfonts,amsbsy,latexsym,dsfont}
\usepackage[textsize=small]{todonotes}       
\usepackage{graphicx}
\usepackage[numeric,initials,nobysame]{amsrefs}
\usepackage{upref,setspace}

\usepackage{caption}
\usepackage{subcaption}
\usepackage{fourier}

\usepackage{enumerate}

\numberwithin{equation}{section}
\numberwithin{figure}{section}
\numberwithin{table}{section}

\newcommand{\cB}{\mathcal{B}}\newcommand{\cC}{\mathcal{C}}

\newcommand{\vzero}{\mathbf{0}}

\newcommand{\ve}{\mathbf{e}}

\newcommand{\vq}{\mathbf{q}}
\newcommand{\vu}{\mathbf{u}}
\newcommand{\vx}{\mathbf{x}}
\newcommand{\vy}{\mathbf{y}}\newcommand{\vz}{\mathbf{z}}

\newcommand{\bR}{\mathbb{R}}

\newcommand{\bZ}{\mathbb{Z}}

\newtheorem{thm}{Theorem}[section]
\newtheorem{lem}[thm]{Lemma}
\newtheorem{cor}[thm]{Corollary}

\newtheorem{ass}[thm]{Assumption}
\newtheorem*{ass*}{Assumption}
\newtheorem*{theorem*}{Theorem}

\theoremstyle{definition}
\newtheorem{rem}{Remark}
\renewcommand{\leq}{\leqslant}
\renewcommand{\geq}{\geqslant}

\newcommand{\eps}{\varepsilon}

\newcommand{\set}[1]{\left\{#1\right\}}

\newcommand{\sss}{\scriptscriptstyle}

\def\qed{ \hfill $\blacksquare$}
\DeclareMathOperator{\E}{\mathbb{E}}
\DeclareMathOperator{\pr}{\mathbb{P}}

\DeclareMathOperator{\var}{Var}

\newcommand{\bpare}[1]{\left( #1 \right)}
\newcommand{\bbrac}[1]{\left[ #1 \right]}

\newcommand{\theline}{(\overline{\vzero, n\ve_1})}
\newcommand{\dlog}{\log \log}
\newtheorem*{lem*}{Lemma}
\DeclareMathOperator{\dist}{d}

\begin{document}

\title[Entropy reduction in Euclidean FPP]{Entropy reduction in Euclidean first-passage percolation}


\author[Damron]{Michael Damron$^1$}
\address{$^1$Georgia Institute of Technology}
\author[Wang]{Xuan Wang$^1$}
\email{mdamron6@gatech.edu, xuanwang9527@gmail.com}

\maketitle

\begin{abstract}
The Euclidean first-passage percolation (FPP) model of Howard and Newman is a rotationally invariant model of FPP which is built on a graph whose vertices are the points of homogeneous Poisson point process. It was shown in \cite{howard2001geodesics} that one has (stretched) exponential concentration of the passage time $T_n$ from $0$ to $n\ve_1$ about its mean on scale $\sqrt{n}$, and this was used to show the bound $\mu n \leq \mathbb{E}T_n \leq \mu n + C\sqrt{n} (\log n)^a$ for $a,C>0$ on the discrepancy between the expected passage time and its deterministic approximation $\mu = \lim_n \frac{\mathbb{E}T_n}{n}$. In this paper, we introduce an inductive entropy reduction technique that gives the stronger upper bound $\mathbb{E}T_n \leq \mu n + C_k\psi(n) \log^{(k)}n$, where $\psi(n)$ is a general scale of concentration and $\log^{(k)}$ is the $k$-th iterate of $\log$. This gives evidence that the inequality $\mathbb{E}T_n - \mu n \leq C\sqrt{\mathrm{Var}~T_n}$ may hold.
\end{abstract}

\section{Introduction}
In \cite{howardnewman1997}, C. D. Howard and C. M. Newman introduced the following Euclidean first-passage percolation (FPP) model on $\bR^d$: Let $Q \subset \bR^d$ be a rate one Poisson point process. Denote by $q(\vx)$, $\vx \in \bR^d$, the closest point to $\vx$ in $Q$, breaking ties arbitrarily. Fix $\alpha > 1$ and define, for $k \geq 1$ and $r = (\vq_1, \cdots, \vq_k)$, a finite sequence of points in $Q$,
\begin{equation*}
	T(r) = \sum_{i=1}^{k-1} \|\vq_i - \vq_{i+1}\|^\alpha,
\end{equation*}
where $\| \cdot \|$ is the Euclidean norm. Such a sequence $r = (\vq_1, \cdots, \vq_k)$ is called a path in $Q$. $r$ can also be viewed as a subset of $Q$ and we write $r \subset Q$.
Define, for $\vq,\vq' \in Q$, $T(\vq,\vq') = \inf\set{T(r)}$, where the infimum is over all finite sequences $r \subset Q$ with $\vq_1 = \vq$ and $\vq_k = \vq'$, and $k$ is the length of $r$. (The condition $\alpha>1$ is imposed because if $0 \leq \alpha \leq 1$, then the straight line segment connecting any two Poisson points is a minimizing path for $T$, and the analysis becomes trivial.) For $\vx, \vy \in \bR^d$, define $T(\vx,\vy) = T(q(\vx),q(\vy))$ and set $T_n = T(0,n\ve_1)$. By subadditivity, the \emph{time constant} $\mu$ exists and is defined by the formula
\[
\mu = \lim_n \frac{\mathbb{E}T_n}{n}.
\]
By the subadditive ergodic theorem, the convergence also holds almost surely, so that in a certain sense, $T_n = \mu n + o(n)$.

In this and related models (lattice FPP and continuum analogues, for example), it is customary to measure the rate of convergence in the definition of $\mu$ by splitting $T_n - \mu n$ into a random fluctuation and nonrandom fluctuation term:
\[
T_n - \mu n = (T_n - \mathbb{E}T_n) + (\mathbb{E}T_n - \mu n).
\]
Typically the random term is analyzed using concentration inequalities (for functions of independent random variables), which lately have developed significantly. In FPP models, current bounds on random fluctuations are still quite far away from the predictions, and this presents an ongoing challenge to researchers. In contrast, there is no general method for providing upper bounds on nonrandom fluctuations of subadditive ergodic sequences. In recent years, though, techniques have been developed \cite{alexander1997, tessera} to bound these nonrandom errors for many lattice models in terms of the random ones. Specifically, if one has a concentration inequality of the type
\begin{equation}\label{eq: taco_head}
\mathbb{P}(|T_n - \mathbb{E}T_n| \geq \lambda \psi(n)) \leq e^{-c\lambda^a}
\end{equation}
for $\lambda \geq 0$ and a suitable function $\psi(n)$ (so far, only results with $\psi(n)$ at least of order $\sqrt{n}$ (in Euclidean FPP) or $\sqrt{n/\log n}$ (in lattice FPP) have been proved), then one can derive the bound
\[
\mu n \leq \mathbb{E}T_n \leq \mu n + C \psi(n) \log n.
\]
(In fact, only the lower tail inequality is usually needed.) A natural question emerges: in these models, can one find $C>0$ such that
\[
\mathbb{E}T_n - \mu n \leq C \sqrt{\mathrm{Var}~T_n}~?
\]
If the answer is yes, it means that the difference $T_n - \mu n$ (used to control geodesics, for instance) can be reasonably well approximated by $T_n - \mathbb{E}T_n$. Furthermore, due to the general lower bounds on nonrandom fluctuations proved in \cite{ADHgamma}, it would suggest that the nonrandom fluctuation term is of the same order as the random one (as is the case in exactly solvable directed last-passage percolation \cite[Corollary~1.3]{baikLPP}).

This question is the focus of our paper. Although we cannot prove this inequality, we show a weaker, but close one. Specifically, our main method is an inductive ``entropy reduction'' technique which shows that for any $k$, there is a constant $C_k$ such that for large $n$,
\[
\mu n \leq \mathbb{E}T_n \leq \mu n + C_k \psi(n) \log^{(k)}n,
\]
where $\psi(n)$ is from \eqref{eq: taco_head} and $\log^{(k)}n$ is the $k$-th iterate of $\log$ (see Theorem~\ref{thm:red-C}). This gives strong evidence that the answer to the above question is yes.

In the next section, we give some background on Euclidean FPP from \cite{howard2001geodesics} and sketch the main strategy to prove general bounds on nonrandom fluctuations in the model. In Section~\ref{sec:main-result}, we state our main assumptions on $\psi$ and the four results (bounds on nonrandom fluctuations, concentration estimates, and geodesic wandering estimates) which come out of our inductive method.

\subsection{Background}
\label{sec:background}

A geodesic between two points $\vx, \vy \in \bR^d$ is a path $r \subset Q$ such that $T(\vx,\vy) = T(r)$. Since $\alpha > 1$, geodesics exist and are unique almost surely \cite[Proposition~1.1]{howard2001geodesics}. Denote by $M(\vx,\vy)$ the geodesic between $\vx$ and $\vy$. Note that $M(\vx, \vy)$ can also be viewed as a subset of $Q$.

First we quote some results from \cite{howard2001geodesics}. Define
\begin{equation}
	\label{eqn:kappa}
	\kappa_1 : = \min\set{1, {d}/{\alpha}}, \mbox{ and } \kappa_2 := {1}/{(4\alpha + 3)}
\end{equation}
and write $\ve_1, \ldots, \ve_d$ for the standard basis vectors of $\mathbb{R}^d$.

\begin{thm}[\cite{howard2001geodesics}, Theorem 2.1]\label{thm:concentration}
	 Define $T_n = T(0, n \ve_1)$. Then there exist constants $C_0, C_1 > 0$ such that $\var(T_n) \leq  C_1 n$ and
	\begin{equation*}
		\pr\bpare{|T_n - \E T_n| > x \sqrt{n}} \leq C_1 \exp( - C_0 x^{\kappa_1}),
	\end{equation*}
	for all $n \geq 0$ and $0 \leq x \leq C_0 n^{\kappa_2}$.
\end{thm}

\begin{thm}[\cite{howard2001geodesics}, Eqn. (4.3)]\label{thm:mean}
	There exists a constant $C_1 > 0$ such that 
	\begin{equation} \label{eqn:mean}
		 n \mu \leq \E T_n \leq n \mu + C_1 \sqrt{n} (\log n)^{1/\kappa_1} .
 	\end{equation}
\end{thm}
Define, for $A,B \subset \bR^d$,
\begin{equation*}
	\dist_{\max}(A,B) = \sup_{\vx \in A} \inf_{\vy \in B} \|\vx-\vy\|.
\end{equation*}
Denote by $\theline$ the line segment between $\vzero$ and $n\ve_1$.
\begin{thm}[\cite{howard2001geodesics}, Theorem 2.4] \label{thm:wandering}
	For any $\eps \in (0, \kappa_2/2)$, there exist constants $C_0, C_1 > 0$ such that 
	\begin{equation*}
		\pr \bpare{\dist_{\max}(M(\vzero, n \ve_1), \theline) > n^{\frac{3}{4}+\eps}  } \leq C_1 \exp(-C_0 n^{3\eps \kappa_1/4}).
	\end{equation*}
\end{thm}
By a simple modification of the proof of \cite[Theorem 2.4]{howard2001geodesics}, one can show that for some constant $C_1 > 0$,
\begin{equation} \label{eqn:wandering}
	\pr \bpare{\dist_{\max}(M(\vzero, n \ve_1), \theline) > C_1 n^{3/4} (\log n)^{1/\kappa_1}  } \to 0 \mbox{ as } n \to \infty.
\end{equation}

The factor $(\log n)^{1/\kappa_1}$ in \eqref{eqn:mean} and \eqref{eqn:wandering} comes from the proof technique. Here we give a sketch of the proof of Theorem \ref{thm:mean}, hinging on the following result, which is \cite[Lemma 4.2]{howard2001geodesics}.
\begin{lem}\label{lem: elusive}
Suppose that the functions $\tau : [0,\infty) \to \mathbb{R}$ and $\sigma: [0,\infty) \to [0,\infty)$ satisfy the following conditions: $\tau(x)/x \to \nu \in \mathbb{R}$, $\sigma(x)/x \to 0$ as $x \to \infty$, $\tau(2x) \geq 2\tau(x) - \sigma(x)$ and $\zeta := \limsup_{x \to \infty} \sigma(2x)/\sigma(x) < 2$. Then for any $c > 1/(2-\zeta)$, $\tau(x) \leq \nu x + c \sigma(x)$ for all large $x$.
\end{lem}
\textbf{Proof:}
The proof is copied from \cite{howard2001geodesics} for completeness. It is easily verified that, for $c> 1/(2-\zeta)$, $\bar{\tau}(x) := \tau(x) - c \sigma(x)$ satisfies $\bar{\tau}(2x) \geq 2 \bar{\tau}(x)$ for all large $x$. Iterating this $n$ times yields $\bar{\tau}(2^n x) \geq 2^n \bar{\tau}(x)$ or $\bar{\tau}(2^n x) / (2^n x) \geq \bar{\tau}(x)/x$. Under our hypotheses on $\tau$ and $\sigma$, $\bar{\tau}(x)/x \to \nu$ as $x \to \infty$, so letting $n \to \infty$ shows that $\bar{\tau}(x)/x \leq \nu$ for all large $x$.
\qed

Returning to the proof of \eqref{eqn:mean}, due to the previous lemma, it suffices to prove $\E T_{2n} \geq 2\E T_n - C_1 \sqrt{n} (\log n)^{1/\kappa_1}$. Now consider the geodesic $M(\vzero, 2n\ve_1)$ and let $\vq$ be the first point in $M(\vzero, 2n\ve_1)$ such that $\|\vq \| \geq n$. Then we have $T_{2n} = T(\vzero, \vq) + T(\vq, 2n\ve_1)$. Then the proof is completed once we show that with positive probability, both of the following bounds hold:
\begin{align*}
	|T_{2n} - \E T_{2n}| &\leq \frac{C_1}{3} \sqrt{n}, \\
	\min\set{T(\vzero, \vq), T(\vq, 2n\ve_1)} &\geq \E T_n -  \frac{C_1}{3} \sqrt{n} (\log n)^{1/\kappa_1}.
\end{align*}
Since $\vq$ is a random point, in order to prove the second bound, one needs to apply Theorem \ref{thm:concentration} to all pairs of the form $(\vzero, \vx)$ where $\vx$ satisfies $\|\vx\| \approx n$. Because we have to apply Theorem \ref{thm:concentration} at least $O(n)$ times, if we use a union bound, we need the probability in Theorem \ref{thm:concentration} to be at most of the order $\frac{1}{n^{r}}$ for some large $r >0$. Taking $x = C_1 (\log n)^{1/\kappa_1}$ in Theorem \ref{thm:concentration} will achieve this and thus complete the sketch of the proof.


Our main goal is to improve the $\log n$ term in Theorem~\ref{thm:mean}. This has been done recently in a lattice FPP model and a directed polymer model in \cite{alexander2011subgaussian,alexander2013subgaussian} by an \textbf{entropy reduction} technique, showing that one can replace the $\log n$ term by $\dlog n$. Their key idea is to exploit the dependence between passage times between nearby points to reduce the number of times a concentration result like Theorem \ref{thm:concentration} is applied. 

The improvement from $\log n$ to $\dlog n$ is important, especially when a sub-gaussian concentration bound for $T_n$ is available. For the lattice FPP model, \cite{damron2015subdiffusive} proved sub-gaussian concentration on the scale of $\sqrt{n/\log n}$ (extending work in \cite{benaimrossignol}). Using this, \cite{alexander2011subgaussian} proved that for a directed FPP model, non-random fluctuations can be bounded by the order
\begin{equation*}
	\sqrt{\frac{n}{\log n}} \cdot \dlog n = o(\sqrt{n}).
\end{equation*}
These bounds have not yet been extended to Euclidean FPP. The strongest concentration inequality to date is Theorem~\ref{thm:concentration} of Howard and Newman.

A consequence of our main results is that one can replace the $\sqrt{n}(\log n)^{1/\kappa_1}$ term in Theorem \ref{thm:mean} to  $\sqrt{n}(\phi(n))^{1/\kappa_1}$ where $\phi(n)$ can be an arbitrary iterate of $\log n$. Our proof works under a general framework which does not depend on any particular scale of concentration. So if a sub-gaussian concentration result for Euclidean FPP is proved, then our result would immediately imply a $o(\sqrt{n})$ bound in Theorem \ref{thm:mean}. \\

\textbf{Notation:} we use bold face letters (e.g. $\vx$, $\vy$, $\vq$) to denote elements in $\bR^d$ or $\bR^{d-1}$. Denote by $\| \cdot \|$ the corresponding $\ell^2$-norm and $\| \cdot \|_\infty$ the $\ell^\infty$-norm. We use $C_0>0$ to denote a small constant and $C_1>0$ a large constant, with values that may vary from case to case. We use notation like $D_{\ref{thm:red-A}}$ to denote constants whose values may depend on $k$ and/or $r$, but not on $n$. The subscript refers to the result number. For example, $D_{\ref{thm:red-A}}$ denotes the constant in Theorem \ref{thm:red-A}.

\section{Main results}
\label{sec:main-result}
In this section, we state the main theorems. We state our results in a general way which does not depend on any one particular concentration result. Let $\psi: (0,\infty) \to (0, \infty)$ be a real function. We assume that we have the following concentration on the scale $\psi(n)$.
\begin{ass}
	\label{ass:concentration}
	There exist constants $C_0>0$, $C_1>0$, $\kappa_1>0$ and $\kappa_2 > 0$ such that
	\begin{equation*}
		\pr \bpare{ |T_n - \E T_n| > \lambda \psi(n)} \leq C_1 \exp\bpare{ - C_0\lambda^{\kappa_1}}
	\end{equation*}
	for all $n \geq 1$ and $ 0 \leq \lambda \leq C_0 n^{\kappa_2}$.
\end{ass}
We put the following assumptions on $\psi$.
\begin{ass}
	\label{ass:phi-psi} 
	There exists $n_0 > 0$ such that $\psi(n)$ is increasing for $n \geq n_0$. In addition, there exist constants $D_{\ref{ass:phi-psi}} > 1$ and $\kappa_3 \in (0,1/2)$ such that for all $n \geq n_0$ and $1 \leq c \leq n^{1/2}$, we have
	\begin{equation*}
		\frac{1}{c^{1-\kappa_3}} \psi(n)  \leq \psi(n/c) \leq \frac{D_{\ref{ass:phi-psi}}}{c^{\kappa_3}} \psi(n).
	\end{equation*}
\end{ass}
Note that the above assumption implies that $\psi(n) = O(n^{1- \kappa_3})$ and $\psi(n)= \Omega(n^\eps)$ for any $\eps \in (0, \kappa_3)$. In addition, the above assumption also implies the following simple bounds: For large $n$ and $ 1 \leq c \leq n^{1/2}$, 
	\begin{equation*}
		\psi(cn) \leq c\psi(n) \mbox{ and } \psi(n/c) \geq \frac{1}{c}\psi(n).
	\end{equation*}

We will assume Assumptions \ref{ass:concentration} and \ref{ass:phi-psi} through out the rest of the paper,
and let constants $C_0$, $C_1$, $\kappa_1$, $\kappa_2$ and $\kappa_3$ be as in Assumptions \ref{ass:concentration} and \ref{ass:phi-psi}.
We further define three constants $\gamma, \beta, \eta > 0$ as follows:
\begin{equation}
	\label{eqn:constant-gamma}
	\gamma := \frac{1}{\kappa_1 \kappa_3}, \quad 
	\beta := \frac{1}{2\kappa_1}, \quad
	\mbox{ and } \eta := \beta + \gamma .	
\end{equation}
These constants show up as exponents in our main theorems below, and reasons for the choices will be clear in the proofs.

Define  $\log^{\sss(0)} n = n$  and  $\log^{\sss (k)} n = \log( \log^{\sss(k-1)} n )$ for $k =1,2, 3,\cdots$, whenever this is well-defined. Write $\vx = (x_1, \vx_2) \in \bR^{d}$ where $x_1 \in \bR$ and $\vx_2 \in \bR^{d-1}$. Define for $n \geq 1$ and $k \geq 0$,
\begin{align*}
	\cB^{\sss (k)}(n) :=& \set{ (x_1, \vx_2) \in \bR^{d}: |x_1| \leq \psi(n), \|\vx_2\| \leq \frac{n^{1/2}\psi^{1/2}(n)}{(\log^{\sss(k)} n )^\eta}}.
\end{align*}



\begin{thm}
	\label{thm:red-A}
	Write $B_1 := \cB^{\sss(k-1)}(n)$ and $B_2 := n\ve_1 + \cB^{\sss(k-1)}(n)$. For any $k \geq 2$ and $r >0$, there exists a constant $D_{\ref{thm:red-A}} =D_{\ref{thm:red-A}}(k,r)  > 0$ such that for large $n$
	\begin{equation*}
		\pr \bpare{ \sup_{\vx,\vx' \in B_1,\; \vy, \vy' \in B_2} |T(\vx,\vy)-T(\vx',\vy')| > D_{\ref{thm:red-A}} \psi(n) }  \leq \frac{1}{(\log^{\sss(k-2)}n)^r}.
	\end{equation*}
\end{thm}

Note that the scale of concentration on Theorem~\ref{thm:red-A} is smaller than that of the next theorem (and is independent of $k$). This is the main reason why we can use estimates for any value of $k$ to give improved ones for $k+1$. 

One key ingredient in the proof of the above result is a simple bound on $|\E T(\vx,\vy) - \E T(\vx',\vy')|$ that reflects the fact that $\E T(\vx, \vy)$ is simply a function of $\|\vx-\vy\|_2$. This is not true for general lattice models. Indeed, it is a standard technique (see \cite{Newman, Kesten}, among many others) to decompose a difference like that from the last theorem as
\begin{align*}
T(\vx,\vy) - T(\vx',\vy') &= [T(\vx,\vy) - \E T(\vx,\vy)] + [\E T(\vx,\vy)-\mu_{\vy-\vx}] \\
&+ [T(\vx',\vy') - \E T(\vx',\vy')] + [\E T(\vx',\vy') - \mu_{\vy'-\vx'}] \\
&+ \mu_{\vy-\vx} - \mu_{\vy'-\vx'}.
\end{align*}
(Here we are writing $\mu_{\vu}$ for the limit $\lim_n \frac{T(0,n\vu)}{n}$, which in our model is simply $\|\vu\|\mu$.) The idea then is to use information about the limiting shape for the model (for instance curvature) to control $\mu_{\vy-\vx}-\mu_{\vy'-\vx'}$ directly, but then one must bound both the random and nonrandom errors on the first two lines. The bounds available for nonrandom errors are generally worse (by some logarithmic factor) than those available for random errors, so one cannot obtain better concentration for $T(\vx,\vy)-T(\vx',\vy')$ than the bounds on nonrandom errors. In our case, we can directly decompose
\[
T(\vx,\vy)-T(\vx',\vy') = [T(\vx,\vy)-\E T(\vx,\vy)] + [T(\vx',\vy') - \E T(\vx',\vy')] + [\E T(\vx,\vy) - \E T(\vx', \vy')],
\]
and exploit the rotational invariance of $\E T$ (from the underlying Poisson process) to obtain bounds without needing control of the nonrandom error.



\begin{thm}
	\label{thm:red-B}
	Write $B_1 := \cB^{\sss(k-1)}(n)$ and $B_2 := n\ve_1 + \cB^{\sss(k-1)}(n)$. For any $k \geq 1$ and $r >0$, there exists a constant $D_{\ref{thm:red-B}} =D_{\ref{thm:red-B}}(k,r)  > 0$ such that for large $n$
	\begin{equation*}
		\pr \bpare{\sup_{\vx \in B_1, \; \vy \in B_2}|T(\vx,\vy) - \E T(\vx,\vy)| > D_{\ref{thm:red-B}} \psi(n) (\log^{\sss(k)} n )^{1/\kappa_1} } \leq \frac{1}{(\log^{\sss(k-1)}n)^r}.
	\end{equation*}
\end{thm}

\begin{thm}
	\label{thm:red-C}
	Let $\mu$ be the time constant. For any $k \geq 1$, there exists a constant $D_{\ref{thm:red-C}} = D_{\ref{thm:red-C}}(k) > 0$ such that for large $n$
	\begin{equation*}
		n \mu \leq \E T(\vzero, n \ve_1) \leq n \mu + D_{\ref{thm:red-C}} \psi(n) (\log^{\sss(k)} (n))^{1/\kappa_1}.
 	\end{equation*} 
\end{thm}
Define for any $\lambda \in \bR$ and $n \geq 1$
\begin{equation*}
	L(\lambda) = L(\lambda, n) := \set{(x_1,\vx_2) \in \bR^d: |x_1 - \lambda| \leq \psi(n)}.
\end{equation*}
Define for $n \geq 1$
	\begin{equation*}
		\bar \cB = \bar \cB(n) := \set{(x_1, \vx_2) \in \bR^d: |x_1| \leq \psi(n), \|\vx_2\| \leq {n^{1/2}\psi^{1/2}(n)}}.
	\end{equation*}
Recall that for $A,B \subset \bR^d$,
\begin{equation*}
	\dist_{\max}(A,B) = \sup_{\vx \in A} \inf_{\vy \in B} \|\vx-\vy\|.
\end{equation*}
\begin{thm}
	\label{thm:red-D}
	Write $\bar B_1 := \bar \cB(n)$ and $\bar B_2 := n \ve_1 + \bar \cB(n)$. For any $k \geq 1$ and $r > 0$, there exists a constant $D_{\ref{thm:red-D}}=D_{\ref{thm:red-D}}(k,r)>0$ such that for all $n$ large and $\lambda \in \bbrac{n/(\log^{\sss(k-1)} n)^{\gamma}, n- n/(\log^{\sss(k-1)} n)^{\gamma}}$
	\begin{equation*}
		\pr \bpare{ \sup_{ \vx \in \bar B_1, \; \vy \in \bar B_2}   \dist_{\max}( L(\lambda) \cap M(\vx,\vy), \theline)  > D_{\ref{thm:red-D}} n^{1/2}\psi^{1/2}(n) (\log^{\sss(k)} n )^\beta } \leq \frac{1}{(\log^{\sss(k-1)}n)^r}.
	\end{equation*}
\end{thm}

We will prove Theorems \ref{thm:red-A} to \ref{thm:red-D} by mathematical induction on $k$. Note that Theorem \ref{thm:red-A} is stated for $k \geq 2$ while the other three theorems are stated for $k \geq 1$. The framework of the mathematical induction can be summarized in the following three steps:
\begin{itemize}
	\item \textbf{Step 1 (Initial): } Prove Theorems \ref{thm:red-B}, \ref{thm:red-C} and \ref{thm:red-D} for $k=1$.
	\item \textbf{Step 2 (Assumption): } Assume that Theorems \ref{thm:red-B}, \ref{thm:red-C} and \ref{thm:red-D} are true for $k = k_0 \geq 1$. Denote these three assumptions by II, III and IV respectively.
	\item \textbf{Step 3 (Induction):} Prove that Theorems \ref{thm:red-A}, \ref{thm:red-B}, \ref{thm:red-C} and \ref{thm:red-D} are true for $k=k_0+1$. Denote these four statements by $\textrm{I}^*$, $\textrm{II}^*$, $\textrm{III}^*$ and $\textrm{IV}^*$ respectively. Then they are proved in the following sequence:
	\begin{align*}
		\textrm{II} + \textrm{III} + \textrm{IV} &\Rightarrow \textrm{I}^*\\
		\textrm{I}^* &\Rightarrow \textrm{II}^* \\
		\textrm{IV} + \textrm{II}^* &\Rightarrow \textrm{III}^* \\
		\textrm{IV} + \textrm{II}^* + \textrm{III}^*  &\Rightarrow  \textrm{IV}^*.
	\end{align*}
\end{itemize}

\textbf{Organization of the paper}: In Section \ref{sec:preliminary-results}, we prove some basic results about the Euclidean FPP model. In Section \ref{sec:initial-step}, we verify the initial step of the mathematical induction. In Section \ref{sec:induction-step}, we complete the induction step of the mathematical induction, and therefore complete the proofs of Theorems \ref{thm:red-A}, \ref{thm:red-B}, \ref{thm:red-C} and \ref{thm:red-D}.

\section{Preliminary Results}
\label{sec:preliminary-results}

In this section, we prove some basic properties about the Euclidean FPP model under the Assumptions \ref{ass:concentration} and \ref{ass:phi-psi}. The proof of these results are analogous to the ones when $\psi(n) = \sqrt{n}$.

As a result of \cite[Lemma 5.2]{howard2001geodesics}, we have the following lemma. Define for $\vx \in \bR^d$ and $n \geq 1$
\begin{equation*}
	B(\vx,n) = \set{\vy \in \bR^d: \|\vy-\vx\|_\infty \leq n}.
\end{equation*}

\begin{lem}\label{lem:qx-x}
	Define the events $F_n$, for $n=1,2,\cdots$, as follows:
	\begin{equation*}
		F_n:=\set{ \forall \vx \in B(\vzero,4n),~ \|\vx-q(\vx)\| \leq \psi^{1/\alpha}(n)}.
	\end{equation*} 
	(i) There exist constants $C_0,C_1 > 0$ such that  
	\begin{equation} \label{eqn:496}
		\pr(F_n^c) \leq C_1 \exp(-C_0 \psi^{d/\alpha}(n)).
	\end{equation}
	(ii) Furthermore, there exists a constant $D_{\ref{lem:qx-x}} > 0$ such that, restricted to $F_n$, we have
	\begin{equation*}
		\sup \set{\|\vq-\vq'\|:  (\vq,\vq') \mbox{ is a geodesic between }  \vq,\vq' \in Q \cap B(\vzero, 4n)} \leq D_{\ref{lem:qx-x}} \psi^{1/\alpha}(n).
	\end{equation*}	
\end{lem}
\textbf{Proof: } (The proof follows exactly from \cite[Lemma 5.2]{howard2001geodesics}, whose statement is similar but with $\psi^{1/\alpha}$ replaced by $n^\gamma$ for some $\gamma \in (0,1)$.) It is sufficient to prove \eqref{eqn:496}. Note that $B(\vzero, 4n)$ can be covered with $O \bpare{\frac{n^d}{\psi^{d/\alpha}(n)}}$ balls of radius $\frac{1}{2}\psi^{1/\alpha}(n)$. If $F_n^c$ occurs, then the intersection of $Q$ and one of these balls must be empty. Therefore
\begin{align*}
	\pr \bpare{F_n^c} 
	\leq& C_1 \cdot \frac{n^d}{\psi^{d/\alpha}(n)} \cdot \exp \bpare{ - 2C_0 \psi^{d/\alpha}(n)} \\
	\leq&  \bbrac{C_1 \cdot \frac{n^d}{\psi^{d/\alpha}(n)} \cdot \exp \bpare{ -  C_0 \psi^{d/\alpha}(n)}} \exp \bpare{ - C_0 \psi^{d/\alpha}(n)}\\
	\leq&  \bbrac{C_1 n^{d-\frac{\kappa_3 d}{2\alpha}} \cdot \exp \bpare{ -  C_0 n^{\frac{\kappa_3 d}{2\alpha}} }} \exp \bpare{ - C_0 \psi^{d/\alpha}(n)},
\end{align*}
where the last line uses the fact that $\psi(n) > n^{\kappa_3/2}$ for large $n$. Then the proof is completed.\qed 

For any $\vx, \vy \in \bR^d$ define $H(\vx,\vy) := \E T(\vx,\vy)$. By the symmetry of the Poisson point process, there is a function $h: \bR^+ \to \bR^+$ such that $H(\vx,\vy) = h(\|\vx-\vy\|)$ where $\|\vx-\vy\|$ is the Euclidean norm. As a result of subadditivity, we have the following simple lemma.

\begin{lem}
	\label{lem:lipschitz}
	There exists a constant $D_{\ref{lem:lipschitz}}> 0$ such that, for all $x,y \geq 0$, $|h(x)-h(y)| \leq D_{\ref{lem:lipschitz}} |x-y|$.
\end{lem}
\textbf{Proof: } By subadditivity, 
\begin{equation*}
	\E T( \vzero, y \ve_1) \leq \E T(\vzero, x \ve_1) + \E T(x \ve_1, y \ve_1).
\end{equation*}
Then since $\E T(x \ve_1, y \ve_1) = \E T(\vzero, (y-x) \ve_1) = h(y-x)$,
\begin{equation*}
	h(y) - h(x) \leq h(y-x).
\end{equation*}
Reversing the roles of $x$ and $y$ gives the same bound for $|h(y)-h(x)|$. Last, we note that an immediate consequence of \cite[Lemma~1]{howardnewman1997} is that $h(x) \leq D_{\ref{lem:lipschitz}} x$ for all $x \geq 0$.
 \qed


We also need the following simple lemma to control the difference of passage times when the endpoints do not differ too much.
\begin{lem}
	\label{lem:compare}
	There exists a constant $D_{\ref{lem:compare}} > 0$ such that, restricted to $F_n$, for $\vx, \vy, \vy' \in B(\vzero, 4n)$ such that $\|\vy - \vy'\| \leq (\psi(n))^{1/\alpha}$,
	\begin{equation*}
		|T(\vx,\vy) - T(\vx, \vy')| \leq D_{\ref{lem:compare}} \psi(n).
	\end{equation*}
\end{lem}
\textbf{Proof: } When restricted to $F_n$, we have $\|q(\vy) - \vy\| \leq \psi(n)^{1/\alpha}$. The proof then follows from the following bound from \cite[(2.14)]{howard2001geodesics}:
	\begin{equation*}
		|T(\vx,\vy) - T(\vx, \vy')| \leq (2\|q(\vy) - \vy\| + 2\|\vy - \vy'\|)^{\alpha}.
	\end{equation*}
\qed

The last result in this section is a global concentration result which plays an important role in verifying the initial cases for the mathematical induction.
\begin{lem}
	\label{lem:concentration}
	Define the set $\cC \subset \bR^d \times \bR^d$ as follows:
	\begin{equation*}
		\cC := \set{ (\vx, \vy) \in \bR^d \times \bR^d: \vx, \vy \in B(\vzero,4n)  \mbox{ and } \|\vx - \vy\| \geq n^{1/2}}.
	\end{equation*}
	For any $r > 0$, there exists a constant $D_{\ref{lem:concentration}} = D_{\ref{lem:concentration}} (r) >0$ such that for all large $n$ 
	\begin{equation*}
		\pr(G_n^c) \leq \frac{1}{n^r},
	\end{equation*}
	where the events $G_n$, $n=1,2,\cdots$ are defined as follows:
	\begin{equation*}
		G_n := \set{ |T(\vx,\vy) - \E T(\vx, \vy)| \leq D_{\ref{lem:concentration}} \psi(n) (\log n)^{1/\kappa_1} \mbox{ for all } (\vx,\vy) \in \cC}.
	\end{equation*}
	
\end{lem}
\textbf{Proof: } For any $(\vx,\vy) \in \cC$, there exists
\begin{equation*}
	(\vx', \vy') \in \cC' := \set{(\vx, \vy) : \vx, \vy \in B(\vzero,4n) \cap \bZ^d  \mbox{ and } \|\vx - \vy\| \geq n^{1/4}}
\end{equation*}
such that $\|\vx-\vx'\| \leq \sqrt{d}$ and $\|\vy-\vy'\|\leq \sqrt{d}$. By Lemma \ref{lem:compare}, restricted to $F_n$, when $n$ is large,
\begin{equation*}
	|T(\vx,\vy) - T(\vx',\vy')| \leq 2D_{\ref{lem:compare}} \psi(n).
\end{equation*}
By Lemma \ref{lem:lipschitz},
\begin{equation*}
	|\E T(\vx, \vy) - \E T(\vx',\vy')| \leq 2D_{\ref{lem:lipschitz}} \sqrt{d}.
\end{equation*}
In the rest of the proof we will replace $D_{\ref{lem:concentration}}$ by $D$ in the definition of $G_n$. Combining the above two bounds, when $n$ is large, $F_n \cap G_n^c$ implies that there exists $(\vx',\vy') \in \cC'$ such that 
\begin{equation*}
	|T(\vx',\vy') - \E T(\vx', \vy')| > D \psi(n) (\log n)^{1/\kappa_1} - 2D_{\ref{lem:compare}} \psi(n) - 2D_{\ref{lem:lipschitz}} \sqrt{d} \geq \frac{D}{2} \psi(n) (\log n)^{1/\kappa_1}.
\end{equation*}
when $D$ is large. By Assumption \ref{ass:concentration}, for any fixed pair $(\vx',\vy') \in \cC'$,
\begin{equation*}
	\pr \bpare{ |T(\vx',\vy') - \E T(\vx', \vy')| > \lambda \psi(\|\vx'-\vy'\|)} \leq C_1 \exp\bpare{-C_0 \min\set{\lambda^{\kappa_1}, (C_0\|\vx' - \vy'\|^{\kappa_2})^{\kappa_1}}}.
\end{equation*}
Let $\lambda = \frac{D \psi(n) (\log n)^{1/\kappa_1}}{2\psi(\|\vx'-\vy'\|)}$. Since $ \|\vx'-\vy'\|_\infty \leq 8 n$ and $n$ is large, one has $\psi(\|\vx'-\vy'\|) \leq 8\psi(n)$ and therefore when $n$ is large,
\begin{equation*}
	\min\set{\lambda, C_0 \|\vx' - \vy'\|^{\kappa_2}} \geq \min\set{ \frac{D}{16} (\log n)^{1/\kappa_1}, {C_0n^{\kappa_2/4}}} = \frac{D}{16} (\log n)^{1/\kappa_1}.
\end{equation*}
Therefore
\begin{equation*}
	\pr \bpare{ |T(\vx',\vy') - \E T(\vx', \vy')| > \frac{D }{2}\psi(n) (\log n)^{1/\kappa_1}}  \leq \frac{C_1}{n^{C_0 (D/16)^{\kappa_1}}}.
\end{equation*}
Since $|\cC'| \leq C_1 n^{2d}$, by a union bound,
\begin{equation*}
	\pr(F_n \cap G_n^c) \leq C_1n^{2d} \cdot \frac{C_1}{n^{C_0 (D/16)^{\kappa_1}}}.
\end{equation*}
Combining this bound with Lemma \ref{lem:qx-x} and taking $D$ large complete the proof. \qed


\section{The Initial Step}
\label{sec:initial-step}

The goal of this section is to verify the initial step of the mathematical induction. Precisely, we will prove the following three lemmas in this section. Lemmas \ref{lem:red-B-initial}, \ref{lem:red-C-initial} and \ref{lem:red-D-initial} imply the $k=1$ cases of Theorems \ref{thm:red-B}, \ref{thm:red-C}, and \ref{thm:red-D} respectively.
Note Lemmas \ref{lem:red-B-initial} and \ref{lem:red-D-initial} are actually stronger than the corresponding initial versions of the theorems.
\begin{lem}
	\label{lem:red-B-initial}
	Define $B_1 = \set{\vx \in \bR^d: \|\vx\|_\infty \leq \psi(n)}$ and $B_2 = n\ve_1 + B_1$. For any $r >0$, there exists a constant $D_{\ref{lem:red-B-initial}} = D_{\ref{lem:red-B-initial}}(r)   > 0$ such that for large $n$
	\begin{equation*}
		\pr \bpare{\sup_{\vx \in B_1, \; \vy \in B_2}|T(\vx,\vy) - \E T(\vx,\vy)| > D_{\ref{lem:red-B-initial}} \psi(n) (\log n )^{1/\kappa_1} } \leq \frac{1}{n^r}.
	\end{equation*}
	In fact, one can take $ D_{\ref{lem:red-B-initial}}(r) = D_{\ref{lem:concentration}}(r)$.
\end{lem}

\textbf{Proof: } When $n$ is large,
\begin{align*}
	&B_1 \subset B(\vzero, 4n) \mbox{ and } B_2 \subset B(\vzero, 4n),\\
	&\vx \in B_1 \mbox{ and } \vy \in B_2 \mbox{ implies } \|\vx - \vy\| \geq n - 2\psi(n) \geq n^{1/2}.
\end{align*}
When $D_{\ref{lem:red-B-initial}}(r) = D_{\ref{lem:concentration}}(r)$, the event considered in this lemma implies $G_n^c$. Therefore Lemma \ref{lem:red-B-initial} follows from Lemma~\ref{lem:concentration} immediately. \qed
\begin{rem}
	Without loss of generality we can assume $\kappa_1$ is so small that $\eta > 1/2$ (recall $\eta$ from \eqref{eqn:constant-gamma}). Then $\frac{n^{1/2}\psi^{1/2}(n)}{n^{\eta}} \leq \psi^{1/2}(n) \leq \psi(n)$, and $\cB^{\sss(0)}(n) \subset  \set{\vx \in \bR^d: \|\vx\|_\infty \leq \psi(n)}$. Therefore Lemma \ref{lem:red-B-initial} implies Theorem \ref{thm:red-B} with $k = 1$.
\end{rem}

\begin{lem}
	\label{lem:red-C-initial}
	There exists a constant $D_{\ref{lem:red-C-initial}} > 0$ such that for large $n$.
	\begin{equation*}
		n \mu \leq \E T(\vzero, n \ve_1) \leq n \mu + D_{\ref{lem:red-C-initial}} \psi(n) (\log n)^{1/\kappa_1}.
 	\end{equation*} 
\end{lem}
\textbf{Proof: } By Lemma~\ref{lem: elusive}, it is sufficient to show that there exists a constant $D >0$ such that for all large $n$
\begin{equation}\label{eqn:723}
	h(2n) \geq 2 h(n) - D  \psi(n) (\log n)^{1/\kappa_1}.
\end{equation}
The proof follows from the proof of \cite[Lemma~4.1]{howard2001geodesics} closely.
Note that restricted to $F_n$, there exists $\vq \in Q \cap \set{\vx \in \bR^d: n- \psi(n) < \|\vx\| < n+ \psi(n)}$ such that $\vq$ is on the geodesic $M(\vzero, 2n\ve_1)$. Therefore
\begin{equation*}
	T(\vzero,\vq) + T(\vq, 2n\ve_1) = T(\vzero, 2n\ve_1).
\end{equation*}
Applying this to an outcome in $F_n \cap G_n$ (which has positive probability), for such a $\vq$ we have
\begin{equation*}
	h(2n) = H(\vzero,2n\ve_1) \geq H(\vzero,\vq) + H(\vq, 2n\ve_1) - 3D_{\ref{lem:concentration}} \psi(n) (\log n)^{1/\kappa_1}.
\end{equation*}
 Then by Lemma \ref{lem:lipschitz}, we have
\begin{equation*}
	\min\set{H(\vzero,\vq), H(\vq, 2n\ve_1)} 
	\geq h(n) - D_{\ref{lem:lipschitz}} \psi(n).
\end{equation*}
Combining the above two inequalities, we have
\begin{equation*}
	h(2n) \geq 2h(n) - 3D_{\ref{lem:concentration}} \psi(n) (\log n)^{1/\kappa_1} -2 D_{\ref{lem:lipschitz}} \psi(n).
\end{equation*}
This implies \eqref{eqn:723} for large $n$. 
\qed

\begin{lem}
	\label{lem:red-D-initial}
	
	Write $\bar B_1 := \bar \cB$ and $\bar B_2 := n \ve_1 + \bar \cB$. For any $r > 0$, there exists a constant $D_{\ref{lem:red-D-initial}}=D_{\ref{lem:red-D-initial}}(r)>0$ such that for all large $n$,
	\begin{equation*}
		\pr \bpare{ \sup_{ \vx \in \bar B_1, \; \vy \in \bar B_2}   \dist_{\max}( M(\vx,\vy), \theline)  > D_{\ref{lem:red-D-initial}} n^{1/2}\psi^{1/2}(n) (\log n )^\beta } \leq \frac{1}{n^r}.
	\end{equation*}
\end{lem}

\textbf{Proof: } Restricted to $F_n$, the event considered in Lemma \ref{lem:red-D-initial} implies that there exist $\vx \in \bar B_1$, $\vy \in \bar B_2$ and $\vq \in Q \cap B(\vzero, 4n)$ such that $ \inf_{\vz \in \theline}\|\vq - \vz\| \geq \frac{D_{\ref{lem:red-D-initial}}}{2} n^{1/2}\psi^{1/2}(n) (\log n)^\beta$ and $\vq$ is on the geodesic from $\vx$ to $\vy$, i.e.,
\begin{equation}
	\label{eqn:758}
	T(\vx,\vq) + T(\vq,\vy) = T(\vx,\vy).
\end{equation}
Meanwhile, elementary geometry shows that there exists a constant $C_0 > 0$ such that for large $n$, $\vx \in \bar B_1$, $\vy \in \bar B_2$ and $\vq \in Q \cap B(\vzero, 4n)$ as above,
\begin{equation*}
	\|\vx-\vq\| + \|\vq-\vy\| -  \|\vx-\vy\| \geq C_0 \frac{(D_{\ref{lem:red-D-initial}} n^{1/2}\psi^{1/2}(n) (\log n )^\beta)^2}{n} =  C_0D_{\ref{lem:red-D-initial}}^2 \psi(n)(\log n)^{1/\kappa_1}.
\end{equation*}
Therefore by Lemma \ref{lem:red-C-initial} and the fact that $\|\vx-\vy\| \leq 2n$,
\begin{align*}
	H(\vx,\vq) + H(\vq,\vy) - H(\vx,\vy) \geq&  \mu\|\vx-\vq\| + \mu\|\vq-\vy\| -  \mu\|\vx-\vy\| - D_{\ref{lem:red-C-initial}} \psi(\|\vx-\vy\|) (\log \|\vx-\vy\|)^{1/\kappa_1} \\
	\geq& \mu\|\vx-\vq\| + \mu\|\vq-\vy\| -  \mu\|\vx-\vy\| - D_{\ref{lem:red-C-initial}} \cdot (2\psi(n))\cdot 2 (\log n)^{1/\kappa_1}\\
	\geq& (C_0D_{\ref{lem:red-D-initial}}^2 -4 D_{\ref{lem:red-C-initial}} )\psi(n)(\log n)^{1/\kappa_1}.
\end{align*}
Comparing \eqref{eqn:758} and the above bound, we have
\begin{equation*}
	\max\set{|T(\vx,\vq) - H(\vx,\vq)|, |T(\vq,\vy) - H(\vq,\vy)|, |T(\vx,\vy) - H(\vx,\vy)|} \geq \frac{C_0D_{\ref{lem:red-D-initial}}^2 -4 D_{\ref{lem:red-C-initial}}}{3} \psi(n)(\log n)^{1/\kappa_1}.
\end{equation*}
Taking $D_{\ref{lem:red-D-initial}}$ so large that $\frac{C_0D_{\ref{lem:red-D-initial}}^2 -4 D_{\ref{lem:red-C-initial}}}{3} > D_{\ref{lem:concentration}}$, the above argument implies that when $n$ is large
\begin{equation*}
	\pr \bpare{ \sup_{ \vx \in \bar B_1, \; \vy \in \bar B_2}   \dist_{\max}( M(\vx,\vy), \theline)  > D_{\ref{lem:red-D-initial}} n^{1/2}\psi^{1/2}(n) (\log n )^\beta } \leq \pr(F_n^c) + \pr(G_n^c),
\end{equation*}
The proof is completed by applying Lemmas \ref{lem:qx-x} and \ref{lem:concentration}. \qed

\begin{rem}
	Theorem \ref{thm:red-D} restricts the geodesic $M(\vzero, n\ve_1)$ to $L(\lambda)$, while Lemma \ref{lem:red-D-initial} removes this restriction in the case $k=1$. Therefore Lemma \ref{lem:red-D-initial} implies Theorem \ref{thm:red-D} with $k=1$.
\end{rem}

\section{The Induction Step}
\label{sec:induction-step}

In this section, we complete the mathematical induction step. We assume that Theorems \ref{thm:red-B}, \ref{thm:red-C} and \ref{thm:red-D} hold for $k = k_0 \geq 1$. Denote these three assumptions by II, III and IV respectively.

	The goal is to prove the $k = k_0 + 1$ cases of Theorems \ref{thm:red-A}, \ref{thm:red-B}, \ref{thm:red-C} and \ref{thm:red-D}. Denote these four statements by $\textrm{I}^*$, $\textrm{II}^*$, $\textrm{III}^*$ and $\textrm{IV}^*$ respectively. Then these four statements are proved in the following sequence:
	\begin{align*}
		\textrm{II} + \textrm{III} + \textrm{IV} &\Rightarrow \textrm{I}^*,\\
		\textrm{I}^* &\Rightarrow \textrm{II}^*, \\
		\textrm{IV} + \textrm{II}^* &\Rightarrow \textrm{III}^*, \\
		\textrm{IV} + \textrm{II}^* + \textrm{III}^*  &\Rightarrow  \textrm{IV}^*.
	\end{align*}
For the ease of reference, we state all assumptions precisely. For simplicity, define $\phi(n) := \log^{\sss(k_0 -1)} n$. Recall the constants $\gamma, \beta, \eta$ from \eqref{eqn:constant-gamma}. Define
\begin{equation*}
	u(n):= \frac{n^{1/2}\psi^{1/2}(n)}{\phi^\eta(n)}, \quad v(n):= n^{1/2}\psi^{1/2}(n) (\log \phi(n))^{\beta}, \quad w(n):= \frac{n}{\phi^{\gamma}(n)},
\end{equation*}
and
\begin{equation*}
	u^*(n):= \frac{n^{1/2}\psi^{1/2}(n)}{(\log \phi(n))^{\eta}}, \quad v^*(n):= n^{1/2}\psi^{1/2}(n) (\dlog \phi(n))^{\beta}, \quad w^*(n):= \frac{n}{(\log \phi(n))^{\gamma}}.
\end{equation*}
Recall that $\bar \cB = \set{(x_1, \vx_2) \in \bR^d: |x_1| \leq \psi(n), \|\vx_2\| \leq n^{1/2} \psi^{1/2}(n)}$. Define
\begin{align*}
	\cB = \cB(n) = \cB^{\sss(k_0-1)}(n) =&  \set{(x_1, \vx_2) \in \bR^d: |x_1| \leq \psi(n), \|\vx_2\| \leq u(n)}, \\
	\cB^* = \cB^*(n) = \cB^{\sss(k_0)}(n) =& \set{(x_1, \vx_2) \in \bR^d: |x_1| \leq \psi(n), \|\vx_2\| \leq u^*(n)}.
\end{align*}

\begin{ass}[$\textrm{II}$]
	\label{ass:ind-b}
	Let $B_1 := \cB$ and $B_2 := n \ve_1 + \cB$.	For any $r >0$, there exists a constant $D_{\ref{ass:ind-b}} =D_{\ref{ass:ind-b}}(r)  > 0$ such that for all large $n$ 
	\begin{equation*}
		\pr \bpare{\sup_{\vx \in B_1, \; \vy \in B_2}|T(\vx,\vy) - \E T(\vx,\vy)| > D_{\ref{ass:ind-b}} \psi(n)(\log \phi(n))^{1/\kappa_1}  } \leq \frac{1}{\phi^r(n)}.
	\end{equation*} 
\end{ass}

\begin{ass}[$\textrm{III}$]
	\label{ass:ind-c}
	Let $\mu$ be the time constant. There exists a constant $D_{\ref{ass:ind-c}} > 0$ such that for large $n$
	\begin{equation*}
		n \mu \leq \E T(\vzero, n \ve_1) \leq n \mu + D_{\ref{ass:ind-c}} \psi(n)(\log \phi(n))^{1/\kappa_1}.
	\end{equation*}
\end{ass}

Recall the definition of $L(\lambda) = L(\lambda,n)$ before Theorem \ref{thm:red-D}.
\begin{ass}[$\textrm{IV}$]
	\label{ass:ind-d}
	Define $\bar{B}_1 = \bar{\mathcal{B}}$ and $\bar{B}_2 = \bar{B}_1 + n\ve_1$. For any $r > 0$, there exists $D_{\ref{ass:ind-d}}=D_{\ref{ass:ind-d}}(r)>0$ such that for large $n$ and $\lambda \in [w(n), n- w(n)]$ we have
	\begin{equation*}
		\pr \bpare{ \sup_{ \vx \in \bar B_1, \; \vy \in \bar B_2} \dist_{\max}( L(\lambda) \cap M(\vx, \vy), \theline)  > D_{\ref{ass:ind-d}} v(n) } \leq \frac{1}{\phi^r(n)}.
	\end{equation*}
\end{ass}

Then we state the four statements that we need to prove in order to complete the mathematical induction as follows.

\begin{lem}[$\textrm{I}^*$]
	\label{lem:ind-a-new}
	Let $B_1^* := \cB^*$ and $B_2^* := n \ve_1 + \cB^*$. For any $r > 0$ there exists a constant $D_{\ref{lem:ind-a-new}} = D_{\ref{lem:ind-a-new}}(r) >0 $ such that for all large $n$
	\begin{equation*}
		\pr \bpare{ \sup_{\vx,\vx' \in B_1^*,\; \vy, \vy' \in B_2^*} |T(\vx,\vy)-T(\vx',\vy')| > D_{\ref{lem:ind-a-new}} \psi(n) }  \leq \frac{1}{\phi^{r}(n)}.
	\end{equation*}
\end{lem}

\begin{lem}[$\textrm{II}^*$]
	\label{lem:ind-b-new}
	For any $r > 0$, there exists a constant $D_{\ref{lem:ind-b-new}} = D_{\ref{lem:ind-b-new}}(r)$	such that for large $n$ 
	\begin{equation*}
		\pr \bpare{\sup_{\vx \in B_1^*, \; \vy \in B_2^*}|T(\vx,\vy) - \E T(\vx,\vy)| > D_{\ref{lem:ind-b-new}} \psi(n) (\log \log \phi(n))^{1/\kappa_1}  } \leq \frac{1}{(\log \phi(n))^r}.
	\end{equation*}
\end{lem}

\begin{lem}[$\textrm{III}^*$]
	\label{lem:ind-c-new}
	There exists a constant $D_{\ref{lem:ind-c-new}} > 0$ such that for large $n$
	\begin{equation*}
		n \mu \leq \E T(\vzero, n\ve_1) \leq n \mu + D_{\ref{lem:ind-c-new}} \psi(n)(\dlog \phi(n))^{1/\kappa_1}.
	\end{equation*}
\end{lem}

\begin{lem}[$\textrm{IV}^*$]
	\label{lem:ind-d-new}
	For any $r > 0$ there exists a constant $D_{\ref{lem:ind-d-new}} = D_{\ref{lem:ind-d-new}}(r) > 0$ such that for large $n$ and $\lambda \in [w^*(n),n - w^*(n)]$ 
	\begin{equation*}
		\pr \bpare{ \sup_{ \vx \in \bar B_1, \; \vy \in \bar B_2}   \dist_{\max}( L(\lambda) \cap M(\vx,\vy), \theline)  > D_{\ref{lem:ind-d-new}}v^*(n) } \leq \frac{1 }{(\log \phi(n))^{r}}.
	\end{equation*}
\end{lem}

One main technique used in the proof is to apply Assumption \ref{ass:ind-b} multiple times, and use many transformed copies of $\cB$ to cover a larger region. More precisely, for any $\vx \in \bR^d$, let $T_{\vx}: \bR^d \to \bR^d$ be the linear transformation such that $T_\vx$ rotates $\ve_1$ to $\frac{1}{\|\vx\|}\vx$ in the plane spanned by $\ve_1$ and $\vx$, and fixes all $\vy$ such that $\vy \perp \vx$ and $\vy \perp \ve_1$. For $\vx \in \bR^d$, define
\begin{equation}
	\theta(\vx) := \arccos \bpare{\frac{\vx \cdot \ve_1}{\|\vx\|}}. \label{eqn:def-theta}
\end{equation}
For any $\vx, \vy \in \bR^d$, define
\begin{equation} \label{eqn:def-tilde-b-x-y}
	\tilde \cB (\vx,\vy) := T_{\vx-\vy} \cB(\|\vx-\vy\|). 
\end{equation}
Note that $\tilde \cB (\vx,\vy)$ is obtained by rotating $\cB(\|\vx-\vy\|)$ by an angle of $\theta(\vx-\vy)$, which maps $\ve_1$ to the direction of $\vx - \vy$. By the symmetry of $\cB(\|\vx-\vy\|)$, we have $T_{\vx-\vy} \cB(\|\vx-\vy\|) = T_{\vy-\vx} \cB(\|\vx-\vy\|)$. Similarly, define
\begin{equation*}
	\tilde \cB^* (\vx,\vy) := T_{\vx-\vy}\cB^*(\|\vx-\vy\|).
\end{equation*}
In the proof of Lemma \ref{lem:ind-a-new}, Assumption \ref{ass:ind-b} is applied to many pairs of boxes of the form $\vx + \tilde \cB (\vx,\vy)$ and $\vy + \tilde \cB (\vx,\vy)$. In the proof of Lemmas \ref{lem:ind-c-new} and \ref{lem:ind-d-new}, 
we also apply Lemma \ref{lem:ind-b-new} to pairs of the form $\vx + \tilde \cB^* (\vx,\vy)$ and $\vy + \tilde \cB^* (\vx,\vy)$. In the rest of this section, we prove some results that control the effect of rotation on such boxes.

\begin{lem}
	\label{lem:rotate}
	Define, for $b \geq a > 0$, $B_{a,b} := \set{(x_1,\vx_2) \in \bR^d: |x_1| \leq a, \|\vx_2\| \leq b}$. Then for $\vz \in \bR^d$,
	\begin{equation*}
		\frac{1}{|\tan \theta(\vz)| \cdot b/a + 1}  B_{a,b} \subset \frac{a}{|\sin \theta(\vz)| \cdot b + |\cos \theta(\vz)| \cdot a}  B_{a,b} \subset T_\vz B_{a,b}.
	\end{equation*}
\end{lem}
Note that the second $``\subset''$ in Lemma \ref{lem:rotate} is optimal, in the sense that the constant $\frac{a}{|\sin \theta(\vz)| \cdot b + |\cos \theta(\vz)| \cdot a}$ can not be improved. The proof of this fact is elementary and therefore omitted. 
Lemma \ref{lem:rotate} immediately implies the following results for $\cB$ and $\cB^*$.

\begin{lem}
	\label{lem:rotate2}
	Suppose $ \sqrt{n} \geq c \geq 1$ and $K > 0$. For any $\vz \in \bR^d$ such that
	\begin{equation}
		\label{eqn:921}
		|\tan \theta(\vz)| \leq \frac{K \psi(n)}{u^*(n)} =  \frac{K \psi^{1/2}(n) (\log \phi(n))^\eta}{n^{1/2}},
	\end{equation}
we have
	\begin{equation*}
		\frac{1}{(1+K) c } \cB(n) \subset T_\vz \cB(n/c), \quad \mbox{ and } \quad \frac{1}{(1+K) c } \cB^*(n) \subset T_{\vz} \cB^*(n/c).
	\end{equation*}
\end{lem}
\textbf{Proof: } First we show that for $ \sqrt{n} \geq c \geq 1$
\begin{equation} \label{eqn:753}
	\frac{1}{c}\cB(n) \subset \cB(n/c) .
\end{equation}
By Assumption \ref{ass:phi-psi} and monotonicity of $\phi(\cdot)$, when $n$ is large,
\begin{align*}
	\psi(n/c) \geq& \frac{1}{c} \psi(n),\\
	u(n/c) = \frac{(n/c)^{1/2}\psi^{1/2}(n/c)} {(\phi(n/c))^\eta} \geq&   \frac{(n/c)^{1/2}(\psi(n)/c)^{1/2}} {(\phi(n))^\eta} = \frac{1}{c}u(n).\\
\end{align*}
This proves \eqref{eqn:753}. Next, by \eqref{eqn:921} and the fact $\log \phi(n) \leq \phi(n)$, we have
\begin{equation*}
	|\tan \theta (\vz)| \leq \frac{K \psi^{1/2}(n) (\phi(n))^\eta}{n^{1/2}} = K \cdot \frac{\psi(n)}{u(n)}.
\end{equation*}
By Lemma \ref{lem:rotate}, this implies that $\frac{1}{K+1}\cB(n) \subset T_\vz \cB(n)$, which combined with \eqref{eqn:753} completes the proof of the statement about $\cB$. The statement about $\cB^*$ can be proved similarly. \qed

Given $\cC \subset \bR^d \times \bR^d$, for $n^{1/2} \geq c \geq 1$ and $K >0$, we say $\cC$ is {\bf $(c, K)$-regular of order $n$}, (or simply $(c, K)$-regular) if for every pair $(\vx,\vy) \in \cC$, we have
\begin{align*}
	&\|\vx-\vy\| \geq n/c, \\
	&|\tan \theta(\vx-\vy)| \leq \frac{K \psi(n)}{u^*(n)} = \frac{K \psi^{1/2} (n)(\log \phi(n))^\eta}{n^{1/2}}.
\end{align*}
Note that in the above definition $c$ and $K$ may also depend on $n$. As a corollary of Lemma \ref{lem:rotate2}: 
\begin{cor}
	\label{cor:rotate}
	If $\cC$ is $(c, K)$-regular of order $n$ for $n^{1/2} \geq c \geq 1$ and $K > 0$, then we have, for every $(\vx,\vy) \in \cC$,
	\begin{equation*}
		\frac{1}{(1+K) c } \cB(n) \subset \tilde \cB(\vx,\vy), \mbox{ and } \frac{1}{(1+K) c } \cB^*(n) \subset \tilde \cB^*(\vx,\vy).
	\end{equation*}
\end{cor}

\textbf{Organization of the rest of this section:} We will prove Lemmas \ref{lem:ind-a-new}, \ref{lem:ind-b-new}, \ref{lem:ind-c-new} and \ref{lem:ind-d-new} in Sections \ref{sec:ind-1}, \ref{sec:ind-2}, \ref{sec:ind-3} and \ref{sec:ind-4} respectively. This will complete the proof of Theorems \ref{thm:red-A}, \ref{thm:red-B}, \ref{thm:red-C} and \ref{thm:red-D}.

\subsection{ $\textrm{II} + \textrm{III} + \textrm{IV} \Rightarrow \textrm{I}^*$}
\label{sec:ind-1}


In this section we prove Lemma \ref{lem:ind-a-new}.

\begin{figure}[!ht]
\centering
\includegraphics[width=6in,trim={0 12cm 0 4cm},clip]{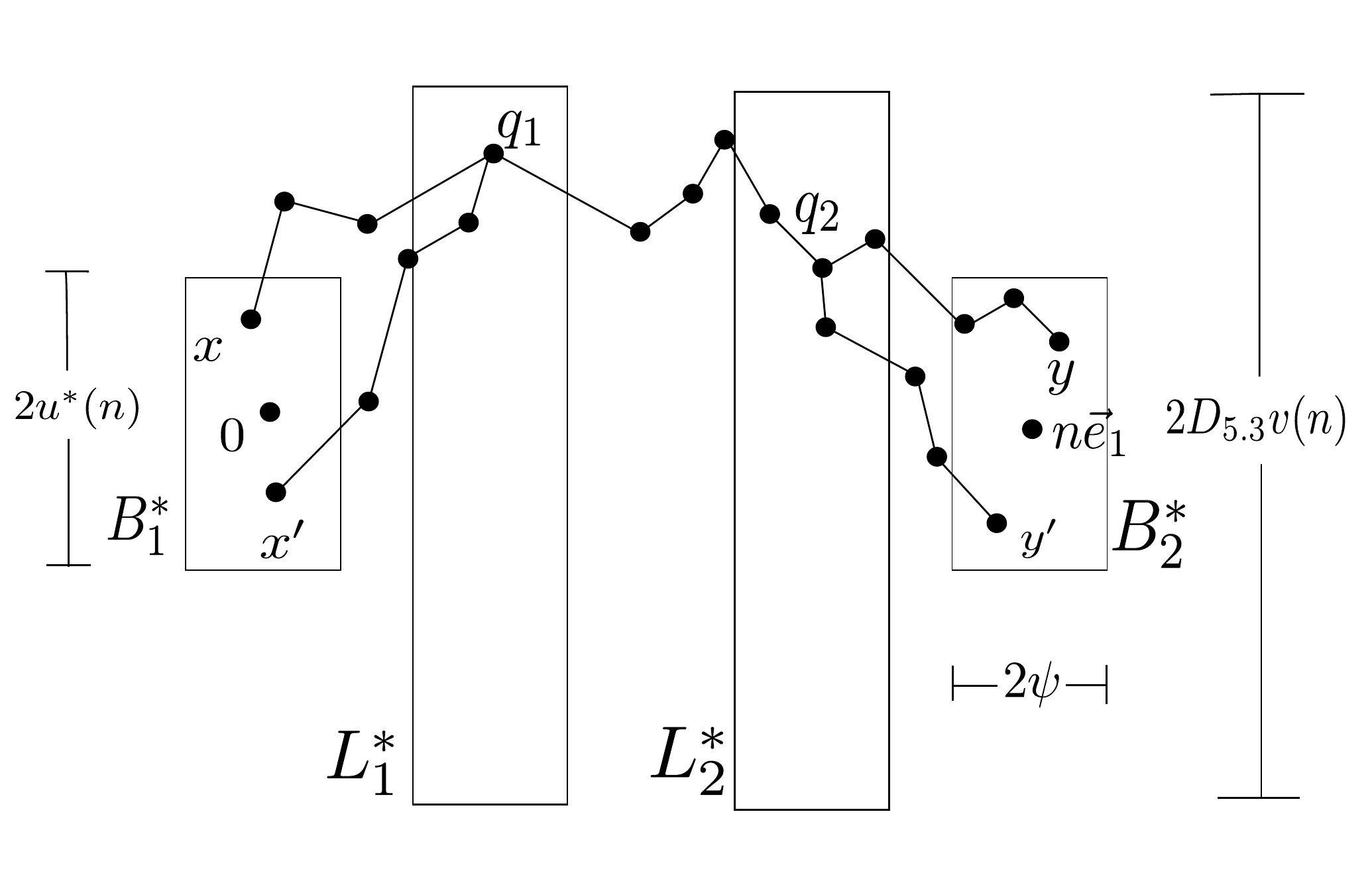}
\caption{Illustration of the proof of Lemma~\ref{lem:ind-a-new}. The path that follows points from $x$ to $q_1$, then to $q_2$, and last to $y$ is a geodesic. One can construct a possibly suboptimal path from $x'$ to $y'$ by taking a geodesic from $x'$ to $q_1$, following the first geodesic from $q_1$ to $q_2$, and then taking a geodesic from $q_2$ to $y'$. Using a similar argument with $x,y$ switched with $x',y'$ produces the main inequality \eqref{eqn:596}.}
\label{fig: fig_1}
\end{figure}

\textbf{Proof of Lemma \ref{lem:ind-a-new}: } Let $r > 0$. Recall that $w^*(n) = n/(\log \phi(n))^\gamma$. Define $L_1 := L(w^*(n))$, $L_2 := L(n - w^*(n))$.  Consider the following events
\begin{equation*}
	H_n := \set{ \sup_{ \vx \in \bar B_1, \; \vy \in \bar B_2} \dist_{\max}( (L_1 \cup L_2) \cap M(\vx, \vy), \theline)  \leq D_{\ref{ass:ind-d}} v(n)  }
\end{equation*}
Since $w^*(n) \in [w(n), n - w(n)]$, by Assumption \ref{ass:ind-d},
we have
\begin{equation}
	\label{eqn:579}
	\pr \bpare{H_n^c} \leq  \frac{1}{\phi^r(n)} +  \frac{1}{\phi^r(n)} 
\end{equation}
Define
\begin{equation}
	L_i^* : = L_i \cap \set{(x_1,\vx_2) \in \bR^d: \|\vx_2\| \leq D_{\ref{ass:ind-d}} v(n) } \mbox{ for } i =1,2.
\end{equation}
Recall $F_n$ from Lemma \ref{lem:qx-x}. When $n$ is so large that $D_{\ref{lem:qx-x}} \psi^{1/\alpha}(n) < 2 \psi(n)$, the event $F_n \cap H_n$ implies that for any $\vx, \vx' \in B_1^* $ and  $\vy, \vy' \in B_2^*$, there exist $\vq_1 \in M(\vx,\vy) \cap L_1^*$ and $\vq_2 \in M(\vx,\vy) \cap L_2^*$ such that
\begin{align*}
	T(\vx,\vy) = T(\vx,\vq_1) + T(\vq_1, \vq_2) + T(\vq_2, \vy), \\
	T(\vx',\vq_1) + T(\vq_1, \vq_2) + T(\vq_2, \vy') \geq T(\vx',\vy').
\end{align*}
Summing up the above two expressions, we have
\begin{equation*}
	T(\vx,\vy) - T(\vx',\vy') \geq [T(\vx,\vq_1) - T(\vx',\vq_1)] + [T(\vq_2, \vy)-T(\vq_2, \vy')].
\end{equation*}
Similarly, there also exists  $\vq_1' \in M(\vx',\vy') \cap  L_1^*$ and $\vq_2' \in M(\vx',\vy') \cap L_2^*$ such that
\begin{equation*}
	T(\vx,\vy) - T(\vx',\vy') \leq [T(\vx,\vq_1') - T(\vx',\vq_1')] + [T(\vq_2', \vy)-T(\vq_2', \vy')].
\end{equation*}
See Figure~\ref{fig: fig_1} for an illustration of the above argument.

 Therefore, restricted to $F_n\cap H_n$,
\begin{equation}
	\label{eqn:596}
	\sup_{\vx,\vx' \in B_1^*,\; \vy, \vy' \in B_2^*} |T(\vx,\vy)-T(\vx',\vy')| \leq \sup_{\vx,\vx' \in B_1^*, \; \vz \in L_1^* }|T(\vx,\vz)-T(\vx',\vz)| + \sup_{\vy, \vy' \in B_2^*, \; \vz \in L_2^* }|T(\vy,\vz)-T(\vy',\vz)|.
\end{equation}
Next, in Lemma \ref{lem:602} we will prove a tail bound for $\sup_{\vx,\vx' \in B_1^*, \; \vz \in L_1^* }|T(\vx,\vz)-T(\vx',\vz)|$. By \eqref{eqn:596}, 
	\begin{align*}
		&\pr \bpare{ \sup_{\vx,\vx' \in B_1^*,\; \vy, \vy' \in B_2^*} |T(\vx,\vy)-T(\vx',\vy')| > 2D_{\ref{lem:602}} \psi(n) }  \\
		\leq& \pr(F_n^c) + \pr(H_n^c) + 2\pr\bpare{\sup_{\vx,\vx' \in B_1^*, \; \vz \in L_1^* }|T(\vx,\vz)-T(\vx',\vz)| > D_{\ref{lem:602}} \psi(n)} \\
		\leq&  C_1\exp( - C_0 \psi^{d/\alpha}(n)) + \frac{2}{\phi^r(n)}  + \frac{2}{\phi^{r-2\eta (d-1) - 1}(n)},
	\end{align*}
	for any fixed $r >0$ and large $n$, where the last line uses Lemma \ref{lem:qx-x}, \eqref{eqn:579} and Lemma \ref{lem:602}. By Assumption \ref{ass:phi-psi}, $\psi(n) = \Omega( n ^{\kappa_3/2})$, so the first two terms in the above display are dominated by the third term. Therefore, for any fixed $r$ and large $n$,
	$$\pr \bpare{ \sup_{\vx,\vx' \in B_1^*,\; \vy, \vy' \in B_2^*} |T(\vx,\vy)-T(\vx',\vy')| > 2D_{\ref{lem:602}} \psi(n) }   \leq \frac{C_1}{\phi^{r-2\eta (d-1) - 1}(n)}. $$
	 Since $r$ can be arbitrarily large, the proof of Lemma \ref{lem:ind-a-new} is completed. \qed

\begin{lem}
	\label{lem:602}
	For any $r > 0$, there exists a constant $D_{\ref{lem:602}} = D_{\ref{lem:602}}(r) > 0$ such that for large $n$,
	\begin{equation*}
		\pr \bpare{\sup_{\vx,\vx' \in B_1^*, \; \vy \in L_1^* }|T(\vx,\vy)-T(\vx', \vy)| > D_{\ref{lem:602}} \psi(n) } \leq \frac{1}{\phi^{r-2\eta (d-1) - 1}(n)}.
	\end{equation*}
\end{lem}
\textbf{Proof: } First, we bound $|\E T(\vx,\vy)-\E T(\vx', \vy)|$. Elementary geometry implies that for large $n$
\begin{align*}
	\sup_{\vx,\vx' \in B_1^*, \vy \in L_1^*}|\|\vx-\vy\|-\|\vx'-\vy\|| 
	\leq& 2 \psi (n) + \frac{D_{\ref{ass:ind-d}} v(n) \cdot 2u^*(n)}{ w^*(n) - 2 \psi(n)} \\
	=& 2 \psi (n) + \frac{D_{\ref{ass:ind-d}} n^{1/2}\psi^{1/2}(n) (\log \phi(n))^\beta \cdot 2n^{1/2}\psi^{1/2}(n) (\log \phi(n))^{-\eta}}{ n (\log \phi(n))^{-\gamma} - 2 \psi(n)} \\
	\leq& 2 \psi(n) + {4 D_{\ref{ass:ind-d}} \psi(n) (\log \phi(n))^{\beta + \gamma - \eta}} \\
	=& (2 + 4D_{\ref{ass:ind-d}}) \psi(n),
\end{align*}
where the third line uses the fact that $\psi(n) = o(n^{1-\kappa_3})$ and so $n (\log \phi (n))^{-\gamma} - 2 \psi (n) > n (\log \phi(n))^{-\gamma}/2$ for large $n$, and the fourth line uses the definition $\eta = \gamma +\beta$. Combining the above bound with Lemma \ref{lem:lipschitz}, we have
\begin{equation}
	\label{eqn:602-1}
	\sup_{\vx,\vx' \in B_1^*, \vy \in L_1^*} |\E T(\vx,\vy) - \E T(\vx', \vy)| \leq (2 + 4D_{\ref{ass:ind-d}}) D_{\ref{lem:lipschitz}} \psi(n).
\end{equation}

Second, we prove the following concentration result: For any $r > 0$ and large $n$,
\begin{equation}
	\label{eqn:602-2}
	\pr\bpare{\sup_{ x' \in B_1^*, y' \in L_1^* }|T(\vx', \vy') - \E T(\vx',\vy')| > 2^{\kappa_3} D_{\ref{ass:phi-psi}}(r) D_{\ref{ass:ind-b}}(r) \psi(n) }
		\leq \frac{1}{\phi^{r-2\eta(d-1) - 1}(n)}.
\end{equation}
To prove this, recall the definition of $\theta(\vx-\vy)$ and $\tilde \cB(\vx,\vy)$ from \eqref{eqn:def-theta} and \eqref{eqn:def-tilde-b-x-y}. 
By Assumption \ref{ass:ind-b}, for every $(\vx,\vy) \in B_1^* \times L_1^*$ and large $n$, 
\begin{equation}
	\label{eqn:615}
	\pr\bpare{\sup_{ \vx' \in \vx + \tilde \cB(\vx,\vy), \vy' \in \vy + \tilde \cB(\vx,\vy) }|T(\vx', \vy') - \E T(\vx',\vy')| > D_{\ref{ass:ind-b}} \psi(\|\vx-\vy\|)[\log \phi(\|\vx-\vy\|)]^{1/\kappa_1} } \leq \frac{1}{\phi^r(\|\vx-\vy\|)}
\end{equation}
When $n$ is large
\begin{align}
	&\|\vx-\vy\| \geq w^*(n) - 2 \psi(n)  \geq \frac{n}{2(\log\phi(n))^{\gamma}} = \frac{w^*(n)}{2}, \label{eqn:837}\\
	&|\tan \theta(\vx-\vy)| \leq \frac{ 2 D_{\ref{ass:ind-d}} v(n)}{w^*(n) - 2 \psi(n) } \leq \frac{4 D_{\ref{ass:ind-d}} \psi^{1/2}(n) \log^{\beta+\gamma}\phi(n)}{n^{1/2}}. \nonumber
\end{align}
Since $\beta + \gamma = \eta$, the set $B_1^* \times L_1^*$ is $(2(\log\phi(n))^{\gamma}, 4 D_{\ref{ass:ind-d}})$-regular. By Corollary \ref{cor:rotate}, we have
\begin{equation}
	\label{eqn:842}
	\cB' := \frac{1}{(1+ 4 D_{\ref{ass:ind-d}}) \cdot 2(\log\phi(n))^{\gamma}} \cB \subset \tilde \cB(\vx,\vy) \mbox{ for all } (\vx,\vy) \in B_1^* \times L_1^*.
\end{equation}

On the other hand, since $\|\vx-\vy\| \leq 2 w^*(n) = \frac{2n}{(\log\phi(n))^{\gamma}} \leq n$ for large $n$, then
\begin{align}
	\psi(\|\vx-\vy\|)[\log \phi(\|\vx-\vy\|)]^{1/\kappa_1} 
	\leq&  D_{\ref{ass:phi-psi}} \bpare{\frac{2}{(\log\phi(n))^{\gamma}}}^{\kappa_3} \psi(n) \cdot (\log\phi(n))^{1/\kappa_1} \nonumber \\
	=& 2^{\kappa_3} D_{\ref{ass:phi-psi}}  {\psi(n)}, \label{eqn:628-1}
\end{align}
where the last line use the relation $\gamma = \frac{1}{\kappa_1 \kappa_3}$.
Since $\phi(\|\vx-\vy\|) \geq \phi(n)/ (2 \log \phi(n))^\gamma$ and $2(\log \phi(n))^{\gamma} \leq \sqrt{n}$ for large $n$,  
\begin{equation}
	\frac{1}{\phi^r(\|\vx-\vy\|)} \leq \frac{2^r(\log \phi(n))^{\gamma r}}{\phi^{r}(n)}. \label{eqn:628-2}
\end{equation}

Combining  \eqref{eqn:842}, \eqref{eqn:628-1} and \eqref{eqn:628-2} in \eqref{eqn:615}, we have for all $(\vx,\vy) \in B_1^* \times L_1^*$,
\begin{equation*}
	\pr\bpare{\sup_{ \vx' \in \vx +  \cB', \vy' \in \vy +  \cB' }|T(\vx', \vy') - \E T(\vx',\vy')| > 2^{\kappa_3} D_{\ref{ass:phi-psi}} D_{\ref{ass:ind-b}} \psi(n) } \leq \frac{2^r(\log \phi(n))^{\gamma r}}{\phi^{r}(n)}.
\end{equation*}
There exists a constant $C_1>0$ such that $B_1^*$ can be covered by 

 $$C_1 \cdot \frac{\psi(n)}{\psi(n)/(\log \phi(n))^{\gamma}} \cdot \bpare{\frac{u^*(n)}{u(n)/(\log \phi(n))^{\gamma}}}^{d-1} = C_1 \phi^{\eta(d-1)}(n) (\log \phi(n))^{\gamma d - \eta(d-1) }$$
copies of $\cB'$, and $L_1^*$ can be covered by 
$$ C_1 \cdot \frac{\psi(n)}{\psi(n)/(\log \phi(n))^{\gamma}} \cdot \bpare{\frac{v(n)}{u(n)/(\log \phi(n))^{\gamma}}}^{d-1} = C_1 \phi^{\eta(d-1)}(n) (\log \phi(n))^{\gamma d + \beta(d-1)}$$ 
copies of $\cB'$. Therefore by a union bound we have
\begin{align*}
	&\pr\bpare{\sup_{ x' \in B_1^*, y' \in L_1^* }|T(\vx', \vy') - \E T(\vx',\vy')| > 2^{\kappa_3} D_{\ref{ass:phi-psi}} D_{\ref{ass:ind-b}} \psi(n) } \\
	\leq& \frac{2^r(\log \phi(n))^{\gamma r}}{\phi^{r}(n)} \cdot  C_1 \phi^{\eta(d-1)}(n) (\log \phi(n))^{\gamma d - \eta(d-1) } \cdot  C_1 \phi^{\eta(d-1)}(n) (\log \phi(n))^{\gamma d + \beta(d-1)}\\
	=& \frac{C_1 (\log \phi(n))^{\gamma r + 2\gamma d - \eta(d-1) + \beta(d-1)}}{\phi^{r-2\eta(d-1)}(n)} \leq \frac{1}{\phi^{r-2\eta(d-1)-1}(n)},
\end{align*}
when $n$ is large. This proves \eqref{eqn:602-2}. Combining \eqref{eqn:602-1} and \eqref{eqn:602-2}, we complete the proof of Lemma \ref{lem:602} with $D_{\ref{lem:602}} = 2^{\kappa_3} D_{\ref{ass:phi-psi}} D_{\ref{ass:ind-b}} +(2 + 4D_{\ref{ass:ind-d}}) D_{\ref{lem:lipschitz}}$. \qed\\


\subsection{$\textrm{I}^* \Rightarrow \textrm{II}^*$}
\label{sec:ind-2}

In this section we prove Lemma \ref{lem:ind-b-new}. Recall that $B_1^* = \cB^*$ and $B_2^* = n\ve_1 + \cB^*$.

\textbf{Proof of Lemma \ref{lem:ind-b-new}: } By the triangle inequality, we have
\begin{align}
	\sup_{\vx \in B_1^*, \; \vy \in B_2^*}|T(\vx,\vy) - \E T(\vx,\vy)| \leq &\sup_{\vx \in B_1^*, \; \vy \in B_2^*}|T(\vx,\vy) - T(\vzero,n\ve_1)| + |T(\vzero,n\ve_1) - \E T(\vzero,n\ve_1)|   \nonumber\\
	&+ \sup_{\vx \in B_1^*, \; \vy \in B_2^*}|\E T(\vx,\vy) - \E T(\vzero,n\ve_1)|.\label{eqn:702}
\end{align}
The first term above can be bounded directly by Lemma \ref{lem:ind-a-new}. The second term can be bounded by the concentration bound in Assumption \ref{ass:concentration}, which implies, for $K = \bpare{{(r+1)}/{C_0}}^{1/\kappa_1}$ and large $n$,
\begin{equation}
	\label{eqn:714-1}
	\pr \bpare{|T(\vzero,n\ve_1) - \E T(\vzero,n\ve_1)| > K \psi(n) (\dlog \phi(n))^{1/\kappa_1} } 
	\leq \frac{C_1}{(\log \phi(n))^{r+1}}.
\end{equation}
To bound the last term in \eqref{eqn:702}, note that for $\vx \in B_1^*$, $\vy \in B_2^*$ and large $n$,
\begin{align*}
	 n - 2 \psi(n) \leq \|\vx-\vy\| \leq& n + 2 \psi(n) + \frac{\bpare{2u^*(n)}^2}{n+2\psi(n)}\\
	 =& n + 2 \psi(n) + \frac{4\bbrac{n^{1/2}\psi^{1/2}(n)/ \log^\eta \phi(n)}^2}{n+2\psi(n)}\\
	 \leq& n + 2 \psi(n) + \frac{4 \psi(n)}{\log^{2\eta}\phi(n)} \leq n + 6\psi(n).
\end{align*}
Then by Lemma \ref{lem:lipschitz}, we have
\begin{equation}
	\label{eqn:714-2}
	\sup_{\vx \in B_1^*, \; \vy \in B_2^*}|\E T(\vx,\vy) - \E T(\vzero,n\ve_1)| \leq 6 D_{\ref{lem:lipschitz}} \psi(n).
\end{equation}
Combining Lemma \ref{lem:ind-a-new}, \eqref{eqn:714-1}, \eqref{eqn:714-2} and \eqref{eqn:702}, when $n$ is large,
\begin{align*}
	& \pr \bpare{\sup_{\vx \in B_1^*, \; \vy \in B_2^*}|T(\vx,\vy) - \E T(\vx,\vy)| > (D_{\ref{lem:ind-a-new}}+K + 6 D_{\ref{lem:lipschitz}}) \psi(n) (\dlog \phi(n))^{1/\kappa_1}  } \\
	\leq& \pr\bpare{\sup_{\vx \in B_1^*, \; \vy \in B_2^*}|T(\vx,\vy) - T(\vzero,n\ve_1)| > D_{\ref{lem:ind-a-new}} \psi(n)} + \pr\bpare{|T(\vzero,n\ve_1) - \E T(\vzero,n\ve_1)| > K \psi(n) (\dlog \phi(n))^{1/\kappa_1}}\\
	\leq& \frac{1}{\phi^r(n)} + \frac{C_1}{\log^{r+1} \phi(n)}  \leq \frac{1}{\log^{r} \phi(n)}.
\end{align*}
The proof of Lemma \ref{lem:ind-b-new} is completed. \qed

\subsection{$\textrm{IV} + \textrm{II}^* \Rightarrow \textrm{III}^*$}
\label{sec:ind-3}

In this section we prove Lemma \ref{lem:ind-c-new}.

\textbf{Proof of Lemma \ref{lem:ind-c-new}: } Write $T_n = T(\vzero, n\ve_1)$ for $n \geq 1$. By Lemma~\ref{lem: elusive}, it suffices to show that there exists a constant $D > 0$ such that for all large $n$,
\begin{equation}
	\label{eqn:771-key}
	\E T_{2n} \geq 2 \E T_n - D \psi(n) (\dlog \phi(n))^{1/\kappa_1}.
\end{equation}
Define for $n \geq 1$,
\begin{equation*}
	L^* := \set{(x_1, \vx_2) \in \bR^d: |x_1 - n | \leq \psi(2n), \; \|\vx_2\| \leq D_{\ref{ass:ind-d}} v(2n)}.
\end{equation*}
For some constant $K >0$ to be decided later, consider the event $E= E_1 \cap E_2 \cap E_3 \cap E_4$ where:
\begin{align*}
	E_1 &= \set{M(\vzero, 2n\ve_1) \cap L^* \neq \emptyset}, \\
	E_2 &= \set{\sup_{\vx \in L^*}|T(\vzero,\vx) - \E T(\vzero,\vx)| \leq K \psi(n) (\dlog \phi(n))^{1/\kappa_1}  }, \\
	E_3 &= \set{\sup_{\vx \in L^*}|T(2n \ve_1,\vx) - \E T(2n \ve_1,\vx)| \leq K \psi(n) (\dlog \phi(n))^{1/\kappa_1},  } \\
	E_4 &= \set{ |T_{2n} - \E T_{2n}| \leq K \psi(n) (\dlog \phi(n))^{1/\kappa_1} }.
\end{align*}
(For the definition of $E_1$, recall that $M(\vzero, 2n\ve_1) \subset Q$.) Restricted to $E_1$, there exists $\vq \in L^*$ such that 
\begin{equation*}
	T(\vzero, 2n\ve_1) = T(\vzero, \vq) + T(\vq, 2n\ve_1).
\end{equation*}
Recall that $\E T(\vx,\vy) = H(\vx,\vy) = h(\|\vx-\vy\|)$. Since $\min\set{\|\vq\|, \|2n\ve_1 - \vq\|} \geq n-\psi(2n) \geq n - 2\psi(n)$, by Lemma \ref{lem:lipschitz}, 
\begin{equation*}
	H(\vzero, \vq) - \E T_n \geq - 2D_{\ref{lem:lipschitz}} \psi(n) \mbox{ and } H(\vq, 2n\ve_1) - \E T_n \geq - 2D_{\ref{lem:lipschitz}} \psi(n).
\end{equation*}
Then restricted to $E$,
\begin{align*}
	\E T_{2n} - 2 \E T_n \geq&  H(\vzero,2n \ve_1) - H(\vzero,\vq) - H(\vq,2n \ve_1) - 4D_{\ref{lem:lipschitz}} \psi(n)\\
	\geq&  T(\vzero, 2n\ve_1) - T(\vzero, \vq) - T(\vq, 2n\ve_1) - 4D_{\ref{lem:lipschitz}} \psi(n) - 3K \psi(n) (\dlog \phi(n))^{1/\kappa_1}\\
	=&- 4D_{\ref{lem:lipschitz}} \psi(n) - 3K \psi(n) (\dlog \phi(n))^{1/\kappa_1}.
\end{align*}
Therefore \eqref{eqn:771-key} holds with $D = 2D_{\ref{lem:lipschitz}} + 3K$ as long as $\dlog \phi(n) > 1$. To complete the proof, it suffices to show that for some choice of $K$, the event $E = \cap_{i=1}^4 E_i$ has positive probability (and therefore is not empty) for all large $n$.

First we  bound $\pr(E_1^c)$. Define $L = \set{(x_1, \vx_2): |x_1 - n| \leq \psi(2n)}$.   By Lemma \ref{lem:qx-x} and Assumption \ref{ass:ind-d} with $n$ replaced by $2n$,
\begin{align*}
	\pr\bpare{ F_n \cap E_1^c } \leq& \pr \bpare{ \dist_{\max}\bpare{L \cap M(\vzero, 2n\ve_1),(\overline{\vzero, 2n\ve_1})} > D_{\ref{ass:ind-d}}v(2n)  } \\
	\leq&  \frac{1}{\phi^{r}(2n)} \to 0 \mbox{ as } n \to \infty.
\end{align*}
Combining this and Lemma \ref{lem:qx-x} we have $\pr(E_1^c) \to 0$ as $n \to \infty$.

Next we bound $\pr(E_4^c)$. By Assumption \ref{ass:concentration}, for $K > 0$ and large $n$
	\begin{align*}
		\pr\bpare{|T_{2n} - \E T_{2n}| > K \psi(n) (\dlog \phi(n))^{1/\kappa_1}} \leq& \pr\bpare{|T_{2n} - \E T_{2n}| > \frac{K}{2} \psi(2n) (\dlog \phi(n))^{1/\kappa_1}} \\
		 \leq& C_1 \exp( - C_0 ({K}/{2})^{\kappa_1} \dlog \phi(n))\\
		 =&  \frac{C_1}{(\log \phi(n))^{C_0 ({K}/2)^{\kappa_1}}} \to 0, \mbox{ as } n \to \infty.
	\end{align*}
 Finally, since $\pr (E_2^c) = \pr(E_3^c)$, we only need to bound $\pr(E_2^c)$. Recall $\tilde \cB^*(\vzero, \vx)$ from \eqref{eqn:def-tilde-b-x-y}. By Lemma \ref{lem:ind-b-new}, for all $\vx \in L^*$
\begin{equation*}
	\pr \bpare{\sup_{\vx' \in \vx + \tilde \cB^*(\vzero, \vx)}|T(\vzero, \vx') - \E T(\vzero, \vx')| > D_{\ref{lem:ind-b-new}} \psi(\|\vx\|) (\dlog \phi(\|\vx\|))^{1/\kappa_1}  } \leq \frac{1}{(\log \phi(\|\vx\|))^r}.
\end{equation*}
Since $n/2 \leq \|\vx\| \leq 2n$ for all $\vx \in L^*$, then the above bound implies for large $n$
\begin{equation}
	\label{eqn:826}
	\pr \bpare{\sup_{\vx' \in \vx + \tilde \cB^*(\vzero, \vx)}|T(\vzero, \vx') - \E T(\vzero, \vx')| > 2 D_{\ref{lem:ind-b-new}} \psi(n) (\dlog \phi(n))^{1/\kappa_1}  } \leq \frac{2^r}{(\log \phi(n))^r}.
\end{equation}

Now we show that the set $\set{\vzero}\times L^*$ is $(2, 8D_{\ref{lem:ind-b-new}})$-regular. Indeed, when $n$ is large, $\|\vx\| \leq 2n$ and
\begin{align}
	v(2n) =&  (2n)^{1/2}\psi^{1/2}(2n) (\log \phi(2n) )^\beta \nonumber \\
	\leq&(2n)^{1/2} \cdot ( 2\psi(n))^{1/2} \cdot 2(\log \phi(n))^\beta  \nonumber\\
	<& 4 v(n). \label{eqn:838}
\end{align}
Then for all $\vx \in L^*$, we have $\|\vx\|\geq n/2$ and
\begin{equation*}
	\tan \theta(\vx) \leq \frac{D_{\ref{ass:ind-d}}v(2n)}{n - \psi(n)} \leq \frac{D_{\ref{ass:ind-d}}\cdot 4 v(n)}{n/2} \leq \frac{8D_{\ref{ass:ind-d}} \psi^{1/2}(\log \phi(n))^\beta}{n^{1/2}}.
\end{equation*}
Thus the set $\set{\vzero}\times L^*$ is $(2, 8D_{\ref{ass:ind-d}})$-regular. Therefore by Corollary \ref{cor:rotate} we have, for all $\vx \in L^*$ and when $D_{\ref{lem:ind-b-new}} > 1$
\begin{equation*}
	\cB' : = \frac{1}{(8D_{\ref{ass:ind-d}} + 1)\cdot 2} \cB^* \subset \tilde \cB^*(\vzero, \vx).
\end{equation*}
Then from \eqref{eqn:826}, we have
\begin{equation}
	\label{eqn:858}
	\pr \bpare{\sup_{\vx' \in \vx + \cB'}|T(\vzero, \vx') - \E T(\vzero, \vx')| > 2 D_{\ref{lem:ind-b-new}} \psi(n) (\dlog \phi(n))^{1/\kappa_1}  } \leq \frac{2^r}{(\log \phi(n))^r}.
\end{equation}
By \eqref{eqn:838}, $L^*$ can be covered by at most 
\begin{equation*}
	C_1 \cdot \frac{\psi(2n)}{\psi(n)} \cdot \bpare{\frac{ v(2n)}{u^*(n)}}^{d-1} \leq C_1 \cdot \frac{2\psi(n)}{\psi(n)} \cdot \bpare{\frac{ 4v(n)}{u^*(n)}}^{d-1} \leq C_1 (\log \phi(n) )^{(\beta+\eta)(d-1)}
\end{equation*}
copies of $\cB'$. Then by the union bound and \eqref{eqn:858}, we have
\begin{align*}
	\pr \bpare{\sup_{\vx \in L^*}|T(\vzero,\vx) - \E T(\vzero,\vx)| > 2 D_{\ref{lem:ind-b-new}} \psi(n) (\dlog \phi(n))^{1/\kappa_1}  } 
	\leq& C_1 (\log \phi(n) )^{(\beta+\eta)(d-1)} \frac{2^r}{(\log \phi(n))^r}. \\
	\leq& \frac{C_1}{(\log \phi(n))^{r-(\beta +\eta)(d-1)}}.
\end{align*}
Taking $K = 2 D_{\ref{lem:ind-b-new}} $ and let $r > (\beta +\eta)(d-1)$, we have $\pr(E_2^c) \to 0$ as $n \to \infty$. Therefore we have proved that $\pr(E^c)$ is small as $n$ is large. The proof of Lemma \ref{lem:ind-c-new} is completed. \qed

\subsection{$ \textrm{IV} + \textrm{II}^* + \textrm{III}^*  \Rightarrow \textrm{IV}^*$}
\label{sec:ind-4}

In this section we prove Lemma \ref{lem:ind-d-new}. 

\begin{figure}[!ht]
\centering
\includegraphics[width=6.5in,trim={0 15cm 2cm 0cm},clip]{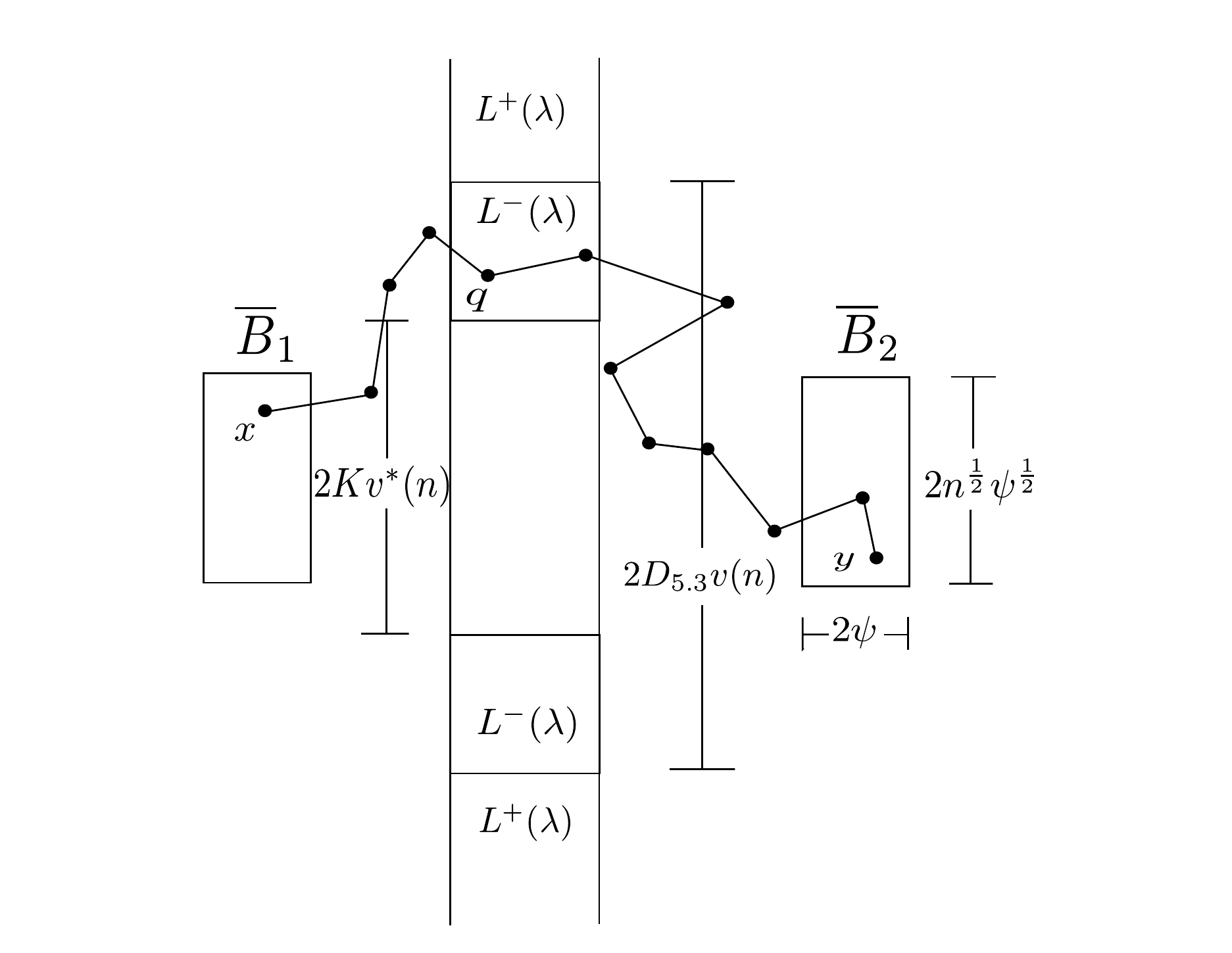}
\caption{Illustration of the event $E(\lambda)$ in the proof of Lemma~\ref{lem:ind-d-new}. The condition is that there are points $\vx \in \overline{B}_1$ and $\vy \in \overline{B}_2$ such that the geodesic $M(\vx,\vy)$ contains a Poisson point in $L^-(\lambda)$. The overall strategy of the proof is to show that geodesics between such points are unlikely to enter $L^+(\lambda)$ (from the event $E^+(\lambda)$) and also unlikely to enter $L^-(\lambda)$ (from the event $E(\lambda)$, illustrated here).}
\label{fig: fig_2}
\end{figure}

\textbf{Proof of Lemma \ref{lem:ind-d-new}: } Let $K$ be a constant whose value will be determined later. Define, for any $\lambda \in [w^*(n),n - w^*(n)]$, 
\begin{align*}
	L^*(\lambda) :=& \set{(x_1,\vx_2) \in \bR^d: |x_1-\lambda| \leq \psi(n), \|\vx_2\| \leq D_{\ref{ass:ind-d}} v(n)},\\
	L^+(\lambda) :=& \set{(x_1,\vx_2) \in \bR^d: |x_1-\lambda| \leq \psi(n), \|\vx_2\| > D_{\ref{ass:ind-d}} v(n)},\\
	L^-(\lambda) :=& \set{(x_1,\vx_2) \in \bR^d: |x_1-\lambda| \leq \psi(n), K v^*(n) < \|\vx_2\| \leq D_{\ref{ass:ind-d}} v(n)},
\end{align*}
Define the events $E^+(\lambda)$ and $E(\lambda)$ for $\lambda \in [w^*(n),n - w^*(n)]$ as follows:
\begin{align*}
	E^+(\lambda) :=& \set{ \exists \vx \in \bar B_1, \vy \in \bar B_2 \mbox{ such that } M(\vx,\vy) \cap L^+(\lambda) \neq \emptyset}.\\
	E(\lambda) :=& \set{ \exists \vx \in \bar B_1, \vy \in \bar B_2 \mbox{ such that } M(\vx,\vy) \cap L^-(\lambda) \neq \emptyset}.
\end{align*}
Then we have
\begin{equation}
	\label{eqn:1295}
	F_n \cap \set{ \sup_{ \vx \in \bar B_1, \; \vy \in \bar B_2}   \dist_{\max}( L(\lambda) \cap M(\vx,\vy), \theline)  > K v^*(n) } \subset E^+(\lambda) \cup E(\lambda).
\end{equation}
By Assumption \ref{ass:ind-d}, for large $n$,
\begin{equation}
	\label{eqn:1043}
	\pr\bpare{ E^+(\lambda)} \leq  \frac{1}{\phi^{r}(n)}. 
\end{equation}
In the rest of the proof, we will prove an upper bound for $\pr(E(\lambda))$. See Figure~\ref{fig: fig_2} for configuration in the event $E(\lambda)$.

Define 
\begin{align*}
	\bar L^- :=  \set{(x_1,\vx_2) \in \bR^d: w^*(n) - \psi(n) \leq x_1 \leq n - w^*(n) + \psi(n), K v^*(n) \leq \|\vx_2\| \leq D_{\ref{ass:ind-d}} v(n)}.
\end{align*}
Further define
\begin{align*}
	\bar \cC :=&  (\bar B_1 \times \bar B_2)  \cup (\bar B_1 \times \bar L^-) \cup (\bar B_2 \times \bar L^-). \\
	\cC(\lambda) :=& (\bar B_1 \times \bar B_2)  \cup (\bar B_1 \times L^-(\lambda)) \cup (\bar B_2 \times L^-(\lambda)), \mbox{ for } \lambda \in \bbrac{w^*(n),n - w^*(n)}
\end{align*}
In order to bound $\pr(E(\lambda))$, we first prove the following relationship: for any $\lambda \in [w^*(n),n - w^*(n)]$,
\begin{equation}
	\label{eqn:1055}
	E(\lambda) \subset E_1(\lambda) := \set{ \sup_{(\vx',\vy') \in \cC(\lambda)} |T(\vx',\vy') - \E T(\vx',\vy')| >  \frac{\mu K^2 - 32D_{\ref{lem:ind-c-new}}}{24} \psi(n) (\dlog \phi(n))^{1/\kappa_1}},
\end{equation}
and then give a bound on $\pr(E_1(\lambda))$ by Lemma \ref{lem:ind-b-new}. Now let us prove \eqref{eqn:1055} first. Note that for any $\vx \in \bar B_1, \vy \in \bar B_2$ and $\vq \in \bar L^-$, elementary geometry shows that when $n$ is large,
\begin{align}
	&\|\vx-\vq\| + \|\vq-\vy\| - \|\vx-\vy\|  \nonumber\\
	\geq& \frac{1}{2} \cdot \frac{\bbrac{K v^*(n) - n^{1/2}\psi^{1/2}(n)}^2}{n/2+\psi(n)}  
	\geq \frac{1}{2} \cdot \frac{\bbrac{K v^*(n) /2}^2}{n} \nonumber\\
	=& \frac{K^2}{8} \psi(n) (\dlog \phi(n))^{2\beta} = \frac{K^2}{8} \psi(n) (\dlog \phi(n))^{1/\kappa_1}. \label{eqn:1287}
\end{align}
By Lemma \ref{lem:ind-c-new} and $\|\vx-\vy\| \leq 2n$,
\begin{align*}
	h(\|\vx-\vy\|) \leq& \mu \|\vx-\vy\| + D_{\ref{lem:ind-c-new}}\psi(\|\vx-\vy\|)(\dlog \phi(\|\vx-\vy\|))^{1/\kappa_1} \\
	\leq& \mu \|\vx-\vy\| + 4 D_{\ref{lem:ind-c-new}}\psi(n)(\dlog \phi(n))^{1/\kappa_1}.
\end{align*}
Combining this and \eqref{eqn:1287} we have
\begin{align}
	h(\|\vx-\vq\|) + h(\|\vq-\vy\|) -h(\|\vx-\vy\|) 
	\geq& \mu \|\vx-\vq\| + \mu \|\vq-\vy\| - \bpare{\mu \|\vx-\vy\| + 4D_{\ref{lem:ind-c-new}} \psi(n)(\dlog \phi(n))^{1/\kappa_1} }  \nonumber\\
	\geq&  \bpare{\frac{\mu K^2}{8} - 4D_{\ref{lem:ind-c-new}}} \psi(n)(\dlog \phi(n))^{1/\kappa_1}. \label{eqn:933}
\end{align}
Since $E(\lambda)$ implies that there exist $\vx \in \bar B_1, \vy \in \bar B_2$ and $\vq \in L^-(\lambda)$ such that 
\begin{equation*}
	T(\vx,\vy) = T(\vx,\vq)+T(\vy,\vq).
\end{equation*}
Combining the above two displays proves \eqref{eqn:1055}.

Next we prove an upper bound for $\pr(E_1(\lambda))$. For any $(\vx,\vy) \in \bar \cC$, define $\tilde B_1 := \vx + \tilde \cB^*(\vx,\vy)$ and $\tilde B_2 := \vy + \tilde \cB^*(\vx,\vy)$. Then by Lemma \ref{lem:ind-b-new},
\begin{equation} \label{eqn:1307}
	\pr \bpare{\sup_{\vx' \in \tilde B_1, \;\vy'  \in  \tilde B_2}|T(\vx',\vy') - \E T(\vx',\vy')| > D_{\ref{lem:ind-b-new}} \psi(\|\vx-\vy\|) (\dlog \phi(\|\vx-\vy\|))^{1/\kappa_1}  } \leq \frac{1}{(\log \phi(\|\vx-\vy\|))^r}.
\end{equation}
Since for large $n$ we have $\|\vx - \vy\| \leq 2n$ and
\begin{align*}
	&\psi(\|\vx-\vy\|) (\dlog \phi(\|\vx-\vy\|))^{1/\kappa_1}  \leq  4  \psi(n) (\dlog \phi(n))^{1/\kappa_1},  \\
	&\frac{1}{(\log \phi(\|\vx-\vy\|))^r} \leq \frac{2^r}{(\log \phi(n))^r}.
\end{align*}
Then \eqref{eqn:1307} implies
\begin{equation} \label{eqn:1316}
	\pr \bpare{\sup_{\vx' \in \tilde B_1, \;\vy'  \in  \tilde B_2}|T(\vx',\vy') - \E T(\vx',\vy')| > 4D_{\ref{lem:ind-b-new}} \psi(n) (\dlog \phi(n))^{1/\kappa_1}  } \leq \frac{2^r}{(\log \phi(n))^r}.
\end{equation}
In addition, when $n$ is large, for all $(\vx, \vy) \in \bar \cC$
\begin{align*}
	& \|\vx-\vy\| \geq w^*(n) - 2\psi(n) \geq \frac{w^*(n)}{2} = \frac{n}{2(\log \phi(n))^\gamma} , \\
	&|\tan \theta(\vx-\vy)| \leq \frac{2 D_{\ref{ass:ind-d}} v(n)}{w^*(n)/2} = \frac{4 D_{\ref{ass:ind-d}} \psi^{1/2}(n) (\log \phi(n) )^{\beta+\gamma}}{n^{1/2}}.
\end{align*}
Since $\beta+\gamma = \eta$, then $\bar \cC$ is $(2(\log \phi(n))^\gamma ,4 D_{\ref{ass:ind-d}} )$-regular. Then by Corollary \ref{cor:rotate}, for all $(\vx,\vy)\in \bar \cC$
\begin{equation*}
	\cB':= \frac{1}{4D_{\ref{ass:ind-d}} + 1} \cdot \frac{1}{2 ( \log \phi(n))^{\gamma}} \cB^* \subset \tilde \cB^*(\vx,\vy).
\end{equation*}
Using this fact in \eqref{eqn:1316} we have
\begin{equation}
	\label{eqn:947}
	\pr \bpare{\sup_{\vx' \in \vx + \cB', \;\vy'  \in \vy + \cB'}|T(\vx',\vy') - \E T(\vx',\vy')| > 4 D_{\ref{lem:ind-b-new}} \psi(n) (\dlog \phi(n))^{1/\kappa_1}  } \leq \frac{2^r}{(\log \phi(n))^r}.
\end{equation}

Note that for $\lambda \in [w^*(n), n - w^*(n)]$,  $\cC(\lambda ) \subset \bar \cC$ and therefore the above bound holds for $(\vx, \vy) \in \cC(\lambda)$. Each of $\bar B_1$ and $\bar B_2$ can be covered by
\begin{equation*}
	C_1 \frac{\psi(n)}{\psi(n)/(\log \phi(n))^{\gamma}} \cdot \bpare{\frac{n^{1/2}\psi^{1/2}(n)}{u^*(n)/(\log \phi(n))^{\gamma}}}^{d-1} = C_1 (\log \phi(n))^{\gamma d + \eta (d-1)}
\end{equation*}
copies of $\cB'$. $L^-(\lambda)$ can be covered by
\begin{equation*}
	C_1 \frac{\psi(n)}{\psi(n)/(\log \phi(n))^{\gamma}} \cdot \bpare{\frac{K v(n)}{u^*(n)/(\log \phi(n))^{\gamma}}}^{d-1} = C_1(\log \phi(n))^{\gamma d + (\eta + \beta) (d-1)}
\end{equation*}
copies of $\cB'$. Then if we take $K$ so large that $(\mu K^2 - 32D_{\ref{lem:ind-c-new}})/24 \geq 4 D_{\ref{lem:ind-b-new}} $ , by \eqref{eqn:947} and the union bound, 
\begin{align*}
	\pr(E_1(\lambda))  
	\leq& \frac{C_1}{(\log \phi(n))^r}  \cdot (\log \phi(n))^{2\gamma d} \bbrac{ (\log \phi)^{(\beta + 2\eta)(d-1)} +  (\log \phi)^{(\beta + 2\eta)(d-1)} + (\log \phi)^{2\eta(d-1)} }\\
	\leq& \frac{C_1}{(\log \phi(n))^{r - 2\gamma d - (2\eta +\beta)(d-1)}}.
\end{align*}
Combine the above bound, \eqref{eqn:1043} and \eqref{eqn:1055} in \eqref{eqn:1295}, taking $K =  \sqrt{(96D_{\ref{lem:ind-b-new}} + 32D_{\ref{lem:ind-c-new}})/\mu}$, we have
\begin{align*}
	&\pr \bpare{ \sup_{ \vx \in \bar B_1, \; \vy \in \bar B_2}   \dist_{\max}( L(\lambda) \cap M(\vx,\vy), \theline)  > K v^*(n) }  \\
	\leq& \frac{1}{\phi^{r}(n)} + C_1 \exp\bpare{-C_0 \psi^{d/\alpha}(n)} +  \frac{C_1}{(\log \phi(n))^{r - 2\gamma d - (2\eta +\beta)(d-1)}}.
\end{align*}
Since $r > 0$ is arbitrary, the proof of Lemma \ref{lem:ind-d-new} is then completed. \qed

\bigskip
\noindent
{\bf Acknowledgements.} The research of M. D. is supported by NSF grant DMS-1419230 and an NSF CAREER grant.

\bibliographystyle{plain}

\end{document}